\documentclass{amsart}
\usepackage{mathptmx}
\usepackage{amssymb}
\usepackage{amsmath}

\newtheorem{theo}{Theorem}[section]

\newtheorem{corol}[theo]{Corollary}
\newtheorem{prop}[theo]{Proposition}
\newtheorem{lem}[theo]{Lemma}
\newtheorem{rem}[theo]{Remark}

\newcommand\RR{{{\mathbb R}}}
\newcommand\cF{{\mathcal F}}

\begin{document}

\title[Analytical regularizing effect
for the Boltzmann equation]
{Analytical regularizing effect
for the radial \\ and spatially homogeneous Boltzmann equation}

\author{LEO GLANGETAS AND MOHAMED NAJEME}

\date{}

\address{Universit\'e de Rouen, UMR-6085, Math\'ematiques
\newline\indent
76801 Saint Etienne du Rouvray, France}

\email{leo.glangetas@univ-rouen.fr}

\address{Universit\'e de Rouen, UMR-6085, Math\'ematiques
\newline\indent
76801 Saint Etienne du Rouvray, France}

\email{mohamed.najeme@univ-rouen.fr}



\begin{abstract}
\noindent In this paper, we consider a class of spatially
homogeneous Boltzmann equation without angular cutoff.
We prove that
any radial symmetric weak solution of the Cauchy problem
become analytic for positive time.
\end{abstract}

\maketitle

\renewcommand{\theequation}{\thesection.\arabic{equation}}
\setcounter{equation}{0}
\section{INTRODUCTION}

This paper deals with the analytic regularity 
of the radially symmetric solutions of the 
following Cauchy problem for the spatially homogeneous Boltzmann equation :
\begin{equation}\label{eq Boltzmann}
\frac{\partial f}{\partial t} = Q(f,f),\quad
v\in \RR^{3},\, t>0; \quad f|_{t= 0} = f_0,
\end{equation}
where
$f(t,v) : \RR^+ \times \RR^3 \longrightarrow \RR$
is the probability density of a gas, 
$v\in\RR^3$ the velocity and
$t \geq 0$ the time.
The Boltzmann collision operator $Q(g,f)$ 
is a bi-linear functional given by 
\begin{equation*}
Q(f,g) = \int_{\RR^3} \int_{S^2}
  B(v-v_*,\sigma) 
\left\{ f(v_*') g(v') - f(v_*) g(v) \right\}d\sigma dv_*,
\end{equation*}
where, for $\sigma \in S^2$,
\begin{equation*}
v'   = \frac{v+v_*}{2} + \frac{|v-v_*|}{2}\sigma, \quad
v'_* = \frac{v+v_*}{2} - \frac{|v-v_*|}{2}\sigma.
\end{equation*}
Theses relations between the post and pre-collisional velocities
follow from the conservation of momentum and kinetic energy.
The non-negative function $B(z,\sigma)$ is called
the Boltzmann collision kernel, depends only on $|z|$
and on the cosine of the deviation angle
$\theta$
\begin{align*}
\cos\theta= \langle\frac{v-v_*}{|v-v_*|},\,\sigma \rangle
\end{align*}
and is defined by
\begin{align*}
B(v-v_*,\cos\theta )
 = \Psi(|v-v_*|) \, b(\cos\theta),\quad
 0 \leq \theta \leq \frac{\pi}{2}.
\end{align*}
We will consider the Maxwellian case 
$\psi \equiv 1$
and we suppose that the cross-section kernel $b$
has a singularity at $\theta = 0$ 
(the so-called non-cutoff problem)
and satisfies :
\begin{equation}\label{Hyp b}
   B(v-v_*,\cos \theta) = b(\cos \theta) 
  \sim |\theta|^{-2-2 s}
\quad\hbox{when}\quad \theta\to0,
\quad 0<s<1.
\end{equation}
We put $\langle v \rangle = \left( 1 + v^2 \right)^{\frac12}$
for $v\in\RR^n$
and we shall use the following standard weighted Sobolev spaces,
for $k, \ell \in \RR$, as
\begin{align*}
L^{p}_{\ell}(\RR^n) &=
\left\{
  f \,\,;\,\, \langle v \rangle^{\ell} f \in L^{p}(\RR^n) 
\right\},
\\
H^k_\ell(\RR^n) &=
\left\{
  f\in {\mathcal S}'(\RR^n) \,\,;\,\, \langle v \rangle^{\ell} f \in H^k(\RR^n)
\right\},
\\
L\log L(\RR^n) &=
\left\{
  f\in {\mathcal S}'(\RR^n) \,\,;\,\,  
 \|f\|_{L\log L} =
\textstyle \int_{\RR^n} |f(v)| \log\left(1+|f(v)|\right) dv<\infty
\right\}.
\end{align*}
The Gevrey space is given for $\alpha>0$ by:
\begin{equation*}
G^{\frac1{\alpha}}(\RR^n) = 
\left\{
  f \,\,;\,\, e^{c_0 \langle D\rangle^\alpha} \in L^2(\RR^n) 
\right\},\\
\end{equation*}
where $ \langle D\rangle = (1 + |D_v|^2 )^{\frac12}$.
Remark that $G^1(\RR^n)$
is the usual analytical functions space.

A solution of Boltzmann equation is known
to satisfy the conservation of mass,
kinetic energy and the entropy inequality:
\begin{align*}
\int_{\RR^3} f(t,v) dv &= \int_{\RR^3} f_0(v) dv, \\
\int_{\RR^3} f(t,v) |v|^2 dv &= \int_{\RR^3} f_0(v) |v|^2 dv, \\
\int_{\RR^3} f(t,v) \log(f(t,v)) dv,
   &\leq \int_{\RR^3} f_0(v) \log(f_0(v)) dv.
\end{align*}
We say that a function $f(v)$ is spatially
radially symmetric
with respect to $v\in\RR^3$
if for any rotation $A$ in $\RR^3$
\begin{equation*}
  f(v) = f(Av).
\end{equation*}

A lot of progress has been made on the study
of the non cut-off problems.  
For the existence of weak solutions, see \cite{V98}
and the references therein.

In \cite{L94}, Lions proved that strong compactness is available
at the level of renormalized solutions.
Then Desvillettes proved in \cite{D95} that
there is a regularizing effect
in the case for radially symmetric solutions 
of the Cauchy problem for a 2D Boltzmann equation
with Maxwellian molecules.
And this is definitively different from the cutoff case,
for which there is no smoothing effect. 
The Sobolev smoothing effect for solutions of the Cauchy problem 
was then studied in other works 
(see \cite{ADVW00, DW04, AMUXY08, HMUY08, MUXY09}).

Some gain of regularity is also obtained
for a solution to the Cauchy problem
of a modified 1D model of the Boltzmann equation 
involving a kinetic transport term (see \cite{DG00}).
For recent works on the non-homogeneous Boltzmann equation,
see \cite{AMUXY11-I, AMUXY11-II, AMUXY11-III}.

In \cite{U84}, Ukai showed that the Cauchy problem
for the Boltzmann equation has a unique
local solution in Gevrey classes.
Then Desvillettes, Furioli and Terraneo proved in \cite{DFT09}
the propagation of Gevrey regularity for solutions of Boltzmann
equation for Maxwellian molecules.
For the non-Maxwellian case,
Morimoto and Ukai considered in \cite{MU10}
the Gevrey regularity of $C^\infty$ solutions 
in the case with a modified kinetic factor 
$\Psi(|v-v_*|) = (1+|v-v_*|^2)^\frac{\gamma}2$
and recently Zhang and Yin in \cite{ZY12}
the case with the general kinetic factor
$\Psi(|v-v_*|)= |v-v_*|^{\gamma}$.
In \cite{MUXY09}, it was proved that 
the solutions of the linearized Cauchy problem
are in the Gevrey space $G^{\frac1{s}}(\RR^3)$
for any $0<s<1$.

Recently, Lekrine and Xu have proved in \cite{LX09}
that,
in the case $0<s<\frac12$,
any symmetric weak solution of the Boltzmann equation belongs
to the Gevrey space $G^{\frac1{2s'}}(\RR^3)$ 
for any $0<s'<s$ and time $t>0$.

In this work, we consider the case $\frac12 \leq s<1$
and we get the following result.
\begin{theo}\label{Th Boltzmann}
Assume that the cross-section kernel $B$ satisfies
\eqref{Hyp b} with $\frac12 < s < 1$ and
the initial datum 
$f_0 \in L^1_{2+2s} \bigcap L\log L(\RR^3)$,
$f_0 \geq 0$ is radially symmetric.
If $f$ is a nonnegative radially symmetric weak solution
of the Cauchy problem 
for the Boltzmann equation \eqref{eq Boltzmann} such that
$f\in L^{\infty}(]0,\infty[; \, L^1_{2+2s} \cap L\log L(\RR^3))$, 
then $f(t,\cdot) \in G^1(\RR^3)$ for any $t>0$.

However, for $s=\frac12$, we have
$f(t,\cdot) \in G^{1/\alpha}(\RR^3)$ for any $0<\alpha<1$ and $t>0$.
\end{theo}

It is well-known that the study of radially symmetric solutions
of the Boltzmann equation can be reduced to the study
of the solutions of the following Kac equation
(see \cite{D95} and also section \ref{section Boltzmann})
\begin{equation}\label{eq Kac}
\left\{\begin{array}{ll}
  & \displaystyle\frac{\partial f}{\partial t}
    = K(f,f) ,\\
  & f|_{t=0} = f_0, \end{array}
\right.
\end{equation}
where $f=f(t,v)$ is the density distribution function
with velocity $v\in\RR$
and the Kac's bilinear collisional operator $K$ is given by
\begin{equation*}
K(f,g) = \int_{\RR} \int_{-\frac{\pi}{2}}^{\frac{\pi}{2}}
  \beta(\theta) 
\left\{ f(v_*') g(v') - f(v_*) g(v) \right\} d\theta dv_*,
\end{equation*}
where
\begin{equation*}
  v' = v \cos\theta - v_* \sin\theta, \quad
  v_*' = v \sin\theta + v_* \cos\theta.
\end{equation*}
The non-negative cross-section $\beta$
satisfies
\begin{equation}\label{beta}
\beta(\theta) = 
  b_0 \frac{|\cos\theta|}{|\sin \theta|^{1+2s}}
\quad\hbox{when }\quad \theta\to 0
\end{equation}
for $0<s<1$ and $b_0>0$.
Remark that
\begin{equation}\label{cond beta}
\int_{-\pi/2}^{\pi/2} \beta(\theta) |\theta|^2 d \theta < \infty.
\end{equation}
There is also conservation of the mass, the kinetic energy and the entropy inequality
for the solutions of the Kac's equation.
We will prove the following result:

\begin{theo}\label{Th Kac}
Assume that the cross-section kernel $\beta$ satisfies
\eqref{beta} with $\frac12 < s < 1$,
the initial datum 
$f_0 \in L^1_{2+2s} \bigcap L\log L(\RR)$.
For $T_0>0$,
if $f\in L^{\infty}([0,\infty[;  L^1_{2+2s}\bigcap L\log L(\RR))$
is a nonnegative weak solution 
of the Cauchy problem of the Kac's equation \eqref{eq Kac},
then $f(t,\cdot) \in G^1(\RR)$ for any $t>0$.

However, for $s=\frac12$, we have
$f(t,\cdot) \in G^{1/\alpha}(\RR)$ for any $0<\alpha<1$ and $t>0$.
\end{theo}

Same as in the paper of \cite{LX09},
the Theorem~\ref{Th Boltzmann} is a direct consequence
of the Theorem~\ref{Th Kac}.
We are reduced to study the Cauchy problem
for spatially homogeneous Kac's equation.

This paper is organized as follows: 
In the next section, we prove some estimates
which will be used in section~\ref{section Kac}.
In section~\ref{section Sobolev},
we study the regularity in weighted Sobovev spaces
for the weak solutions of the Cauchy problem of the Kac's equation.
The section~\ref{section Kac} is devoted to the proof
of the Theorem \ref{Th Kac}
and in section~\ref{section Boltzmann} we 
conclude the proof the Theorem~\ref{Th Boltzmann}.

\renewcommand{\theequation}{\thesection.\arabic{equation}}
\setcounter{equation}{0}
\section{ESTIMATES OF THE COMMUTATORS}\label{section Estimate}
In this section, we will get the estimates
of some terms that we call ``commutators'' and
we will see in section \ref{section Kac} that 
they are the main point
to get the regularity of weak solutions for the Cauchy problem
of the Kac's equation. 
We recall the following coercivity inequality
deduced from the non cut-off of collision kernel.
\begin{prop}\label{Prop coer} (see \cite{ADVW00})
Assume that the cross-section
and satisfies the assumption \eqref{beta}.
Let $f \geq 0$, $f\neq 0$, 
$f\in L^1_1(\RR) \bigcap L\log L(\RR)$, 
then there exists a constant $c_f>0$, 
depending only of, 
$\beta$, $\|f\|_{L^1_1}$ and $\|f\|_{L\log L}$, such that 
\begin{equation*}
-\left(K(f,g),\,\,g \right)_{L^2}
\geq 
c_f\|g\|^2_{H^{s}}\,\, 
  - \,\,C\|f\|_{L^1}\|g\|^2_{L^2} 
\end{equation*}
for any smooth function $g\in H^s(\RR)$.
\end{prop}
\noindent
{\bf Remark.} 
From \cite{HMUY08, MUXY09}, 
if $m,\ell\in \RR$, $0<s<1$ and $f$ and $g$ are suitable functions, 
the Kac collision kernel has the following regularity ($\ell^+ = \max(0, \ell)$)
\begin{equation*}
\| K(f,g) \|_{H^m_\ell}
\leq
C
\| f \|_{L^1_{\ell^+ + 2 s}} \,
\| g \|_{H^{m+2 s}_{(\ell+2s)^+}}. 
\end{equation*}
As in \cite{MUXY09}, we introduce the following mollifier
\begin{equation}\label{G_delta}
   G_\delta(t, \xi) =
  \frac{     e^{c_0 t \langle \xi\rangle^{\alpha} } }
  {1+ \delta e^{c_0 t \langle \xi\rangle^{\alpha} } }
\end{equation}
where
$\langle \xi\rangle = (1+|\xi|^2)^{\frac12}$,
$\xi\in \RR$, $c_0>0$ and $0<\delta<1$
will be chosen small enough
and $\alpha\in]0,2[$ are fixed.
It is easy to check that, for any $0<\delta<1$,
\begin{equation*}
G_{\delta}(t,\xi) \in L^{\infty}(]0,T[\times \RR).
\end{equation*}
We denote by $\hat{f}$ the Fourier transform of $f$
\begin{equation*}
\hat{f} (\xi) 
= \cF (f) (\xi)
= \int_{\RR}e^{-i v.\xi} f(v) \, dv
\end{equation*}
and by $G_{\delta}(t,D_v)$
the Fourier multiplier of symbol $G_{\delta}(t,\xi)$
(see \cite{H85})
\begin{equation*}
G_{\delta} \, g(t,v) 
  = G_{\delta}(t,D_v) g(t,v) 
  = \cF^{-1}\left( G_{\delta}(t,\cdot)\hat{g}(t,\cdot)\right)(v).
\end{equation*}
The proof of Theorem \ref{Th Kac}
will be based on the uniform estimate
with respect to $0<\delta<1$ of
$\|G_{\delta}(t,D_v)f(t,\cdot)\|_{L^2_1}$
where $f(t,\cdot)$ is a weak solution of the
Cauchy problem of the Kac's equation \eqref{eq Kac}.

In the following, $C$ will represent a generic constant independent of $\delta$
and $t\in[0,T]$ (but it will depend on the kernel $\beta$ and the norms
$\|f(t,\cdot)\|_{L^1_2}$, $\|f(t,\cdot)\|_{L \log L}$ used for the coercivity).
\begin{lem}\label{lem estim1 G}
Let $T>0$. 
We have that for any $0<\delta<1$ 
and $\leq t\leq T, \xi \in \RR$,
\begin{align*}
\left| \partial_t G_\delta(t, \xi) \right| 
&\leq 
  c_0 \langle \xi\rangle^{\alpha} G_\delta(t, \xi),
\\
\left| \partial_\xi G_\delta(t, \xi) \right| 
&\leq 
 \alpha c_0 t \langle \xi\rangle^{\alpha-1} G_\delta(t, \xi),
\\
\left| \partial ^2_{\xi} G_\delta(t, \xi) \right| 
&\leq 
 C \langle \xi\rangle^{2 \alpha - 2} G_\delta(t, \xi)
\end{align*}
with $C>0$ independent of $\delta$ and $t\in[0, T]$.
\end{lem}
\begin{proof}
We compute
\begin{align*}
\partial_{t} & G_{\delta}(t,\xi) 
= 
c_0 \frac{\langle \xi\rangle^{\alpha}}
  {1+ \delta e^{c_0 t \langle \xi\rangle^{\alpha}}} 
  \, G_{\delta}(t,\xi),
\\
\partial_{\xi} & G_{\delta}(t,\xi) 
= \alpha c_0 t \, \xi \, (1+|\xi|^2)^{\frac{\alpha}2-1}
  G_{\delta}(t,\xi)
  \frac1{1+\delta e^{c_0 t \langle \xi\rangle^{\alpha}}},
\\
\partial^2_{\xi} & G_\delta(t, \xi)
  = 
  \left( 
  \alpha c_0 t \, \xi \,
  (1+|\xi|^2)^{\frac{\alpha}2-1} \right)^2
  G_{\delta}(t,\xi)
  \frac{1-\delta e^{c_0 t \langle \xi\rangle^{\alpha}}}
       {\left(1+\delta e^{c_0 t \langle \xi\rangle^{\alpha}}\right)^2}
\\
  &+
  \alpha c_0 t 
  \left( 
   (1+|\xi|^2)^{\frac{\alpha}2-1} 
 + (\alpha - 2) \, \xi^2 \, 
  (1+|\xi|^2)^{\frac{\alpha}2-2} 
  \right)
  G_{\delta}(t,\xi)
  \frac1{1+\delta e^{c_0 t \langle \xi\rangle^{\alpha}}}, \nonumber
\end{align*}
and the estimates of the lemma follow easily.
\end{proof} 
\begin{lem}\label{lem estim2 G}
There exists $C>0$ such that for all
$0<\delta<1$ and $\xi \in \RR$
\begin{align*}
\left| G_\delta(t, \xi) - G_\delta(t, \xi \cos \theta) \right| 
  &\leq 
   C \sin^2 \textstyle{\frac{\theta}2} \langle \xi\rangle^{\alpha} 
   G_\delta(t, \xi \cos \theta) G_\delta(t, \xi \sin \theta),
\\
  \left|  (\partial_{\xi} G_\delta)(t, \xi) 
        - (\partial_{\xi} G_\delta)(t, \xi \cos \theta) \right| 
  &\leq 
    C \sin^2 \textstyle{\frac{\theta}2} \langle \xi \rangle^{2 \alpha - 1} G_\delta(t, \xi \cos \theta) G_\delta(t, \xi \sin \theta).
\end{align*}
\end{lem}
\begin{proof}
This lemma \ref{lem estim2 G} is proved by Taylor formula,
the estimates from lemma \ref{lem estim1 G} and the
following inequality :
\begin{equation}\label{G}
G_{\delta}(t,\xi) \leq 
  3 G_{\delta}(t,\xi \cos \theta) G_{\delta}(t,\xi \sin\theta).
\end{equation}
\end{proof} 
We now estimate the commutator of the Kac's operator with the mollifier:
\begin{prop}\label{Prop estim com0}
Assume that $0<\alpha<2$.
Let $f,g \in L^2_1$ and $h\in H^{\alpha/2}(\RR)$,
then we have
\begin{equation*}
| \big(G_\delta K(f, g) , h\big)_{L^2} 
- \big(K(f, G_\delta g) , h\big)_{L^2}  |
\leq
  C \|G_\delta f\|_{L^2_1}
  \|G_\delta g\|_{H^{\alpha/2}} \|h\|_{H^{\alpha/2}}.
\end{equation*}
\end{prop}
\begin{proof}
By definition, of $G_\delta$ we have for a regular $f$,
\begin{equation*}
\cF(G_{\delta}f)(\xi) = G_{\delta}\hat{f}(\xi),
\end{equation*}
and
\begin{equation*}
\cF(vG_{\delta}f)(\xi) =
 i\partial_{\xi}\left( G_{\delta}(t,\xi)\hat{f}(t,\xi)\right).
\end{equation*}
We recall the Bobylev formula
\begin{equation}\label{bobylev}
\cF\left(K(f,g) \right) (\xi)
= 
\int_{-\frac{\pi}{2}}^{\frac{\pi}{2}}
\beta(\theta)
\left\lbrace 
  \hat{f}(\xi \sin \theta)\hat{g}(\xi\cos\theta)-\hat{f}(0)\hat{g}(\xi)
\right\rbrace 
d\theta.
\end{equation}
From the Bobylev and Plancherel formulas
\begin{align*}
\bigl( &G_\delta K(f, g) , h \bigr)_{L^2}
- 
\bigl(  K(f, G_\delta g) , h \bigr)_{L^2}
\\
&= 
\int_{\RR_{\xi}} \int_{-\frac{\pi}{2}}^{\frac{\pi}{2}}
  \beta(\theta) G_\delta(t,\xi) 
  \left\lbrace 
    \hat{f}(\xi \sin \theta)\hat{g} (\xi\cos\theta) -\hat{f}(0) \hat{g}(\xi) 
  \right\rbrace 
  d\theta \overline{\hat{h}} d\xi
\\
&-
  \int_{\RR_{\xi}}\int_{-\frac{\pi}{2}}^{\frac{\pi}{2}}\beta(\theta)
  \left\lbrace
    \hat{f}(\xi \sin \theta) \, \cF(G_\delta g)(\xi \cos \theta)
  - \hat{f}(0) \, \cF(G_\delta g)(\xi ) 
  \right\rbrace 
  \overline{\hat{h}(\xi)}d\theta d\xi 
\\
&= \int_{\RR_{\xi}} \int_{-\frac{\pi}{2}}^{\frac{\pi}{2}}
  \beta(\theta)\hat{f}(\xi \sin \theta)
  \left\lbrace G_{\delta}(\xi) - G_{\delta}(\xi \cos \theta) \right\rbrace 
  \hat{g}(\xi\cos\theta)\overline{\hat{h}(\xi)} d\theta d\xi.
\end{align*}
By the previous formula, lemma \ref{lem estim2 G}
and the Cauchy-Schwarz inequality we have 
\begin{align*}
\Bigl|\bigl( G_\delta K(f, g) , h\bigr)_{L^2}
    - \big(K(f, & G_\delta g) , h\big)_{L^2}\Bigr|
\\
\leq 
  & \int_{\RR_{\xi}}
\int_{-\frac{\pi}{2}}^{\frac{\pi}{2}}
\beta(\theta)\sin^2\textstyle{\frac{\theta}2}
|G_{\delta}(\xi\sin\theta)\hat{f}(\xi\sin\theta)|
\\
  & \times
|G_{\delta}(\xi\cos\theta)\hat{g}(\xi\cos\theta)|\langle \xi \rangle^{\alpha} | \hat{h}(\xi) | d\theta d\xi
\\
\leq & C \|G_{\delta} \hat f\|_{L^{\infty}} 
  \|G_{\delta} g\|_{H^{\alpha/2}}
  \|h\|_{H^{\alpha/2}}
\\
\leq & C \, \|G_{\delta} f \|_{L^2_1}
 \|G_{\delta} g\|_{H^{\alpha/2}} \|h\|_{H^{\alpha/2}}
\end{align*}
where we have used the following continuous embedding
\begin{equation*}L^2_1(\RR) \subset L^1(\RR)\end{equation*}
and the assumption \eqref{cond beta} on the kernel $\beta$.
\end{proof} 
We again estimate the commutator of the Kac's operator
with the mollifier weighted as in \cite{LX09}.
We will need to use a property of symmetry for the Kac's operator.
\begin{prop}\label{Prop estim com1}
Assume that $\frac12 < s<1$ and let
$f, g \in L^1_2(\RR)$ and $h\in H^{\frac12}(\RR)$.
Then we have
\begin{align*}
\big|\big( (v G_\delta) K(f, g) , h \big)_{L^2} 
   &- \big( K(f, (v G_\delta) g)  , h \big)_{L^2} \big|
\\
&\leq\,
C \left(
\|f\|_{L^1_2} + C \|G_\delta f\|_{L^2_1}
\right)
\,
\|G_\delta g \|_{H^{\frac12}_1}
\|h\|_{H^{\frac12}}.
\end{align*}
\end{prop}
\noindent
{\bf Remark}.
For $s=\frac12$, the previous estimate is not enough accurate.
In order to use some interpolation argument, we will need the following
estimate.
\begin{prop}\label{Prop estim com1 s=1/2}
Assume that $s=\frac12$ and let $0<\alpha, \alpha'<1$,
$f, g \in L^1_2(\RR)$,
and $h\in H^{\frac{\alpha}2}(\RR)$.
Then we have
\begin{align*}
 \big| \big( (v G_\delta) &K(f, g) , h \big)_{L^2}
 - \bigl(K(f, (v G_\delta) g) , h \bigr)_{L^2} \big|
 \\
 &\leq\,
C \|f\|_{L^1_{1+{\alpha'}}}
\|G_\delta g \|_{H^{\frac{{\alpha'}}2}}
  \|h\|_{H^{\frac{{\alpha'}}2}}
+ C \|G_\delta f\|_{L^2_1}
\|G_\delta g \|_{H^{\frac{\alpha}2}}
\|h\|_{H^{\frac{\alpha}2}}.
\end{align*}
\end{prop}

We will prove these Propositions by using
the Bobylev formula \eqref{bobylev}
and the Plancherel formula.
We can write
\begin{align*}
\big( (v G_\delta) K(f, g) , h\big)_{L^2}
 - \big(K(f, (v G_\delta) g) , h\big)_{L^2} 
= 
i \int_{\RR_{\xi}}\int_{-\frac{\pi}{2}}^{\frac{\pi}{2}}
\beta(\theta) 
\, A(\xi,\theta) \,
\overline {{\hat h}(\xi)} \,d\theta d\xi 
\end{align*}
where
\begin{equation*}
A(\xi,\theta) = 
\partial_\xi 
\left\{
  \hat f(\xi\sin\theta) G_\delta(\xi) \hat g(\xi\cos\theta)
\right\}
-
\hat f(\xi\sin\theta)
\partial_\xi 
\left\{
  G_\delta \, \hat g
\right\}(\xi\cos\theta).
\end{equation*}
We decompose $A = A_1 + A_2 + A_3$
where
\begin{align*}
&A_1 = \sin\theta \,(\partial_\xi \hat f)(\xi \sin\theta) 
        \, G_\delta(\xi) \, \hat g(\xi\cos\theta),
\\
&A_2 = \hat f(\xi \sin\theta)
  \left\{
  G_\delta(\xi) \cos\theta 
 -G_\delta(\xi\cos\theta) 
  \right\} \, (\partial_\xi\hat g)(\xi\cos\theta),
\\
&A_3 = \hat f(\xi \sin\theta)
  \left\{
  \partial_\xi G_\delta (\xi) 
 -(\partial_\xi G_\delta) (\xi\cos\theta)
  \right\}
\, \hat g(\xi\cos\theta),
\end{align*}
and we put for $k=1, 2, 3$
\begin{equation*}
I_k = i \int_{\RR_{\xi}}\int_{-\frac{\pi}{2}}^{\frac{\pi}{2}} \beta(\theta)
\, A_k(\xi,\theta) \,
\overline {{\hat h}(\xi)} \,d\theta d\xi.
\end{equation*}
Therefore we have
\begin{equation}\label{com1<=}
\big| \big( (v G_\delta) K(f, g) , h\big)_{L^2} 
  - \big(K(f, (v G_\delta) g) , h\big)_{L^2} \big|
\leq |I_1| + |I_2| + |I_3|.
\end{equation}
In the following, we will estimate the three
terms $I_1$, $I_2$ and $I_3$.

\medskip
{\bf Estimate of $I_1$.}
\medskip
We decompose $I_1 = I_{1a} + I_{1b}$ where
\begin{align*}
&I_{1a} = 
i \int_{\RR_{\xi}}\int_{-\frac{\pi}{2}}^{\frac{\pi}{2}}
\beta(\theta) \, \sin\theta \, (
\partial_\xi \hat f)(\xi \sin\theta) 
G_\delta(\xi \cos\theta) \hat g(\xi \cos\theta)
\overline {{\hat h}(\xi)} d\theta d\xi,
\\ 
&I_{1b} = 
 i\int_{\RR_{\xi}}\int_{-\frac{\pi}{2}}^{\frac{\pi}{2}} \beta(\theta) \, \sin\theta \, (\partial_\xi \hat f)(\xi \sin\theta)
\left(G_\delta(\xi ) - G_\delta(\xi\cos\theta ) \right)
\hat g(\xi \cos\theta)
\overline {{\hat h}(\xi)} d\theta d\xi.
 \end{align*} 
\begin{lem}\label{lem I1a}
Suppose that $\frac12 < s < 1$. Then there exists a constant $C$
such that 
\begin{equation*}
|I_{1a}| \leq
C \|f\|_{L^1_2}
\| G_\delta g \|_{H^{\frac12}}
\| h \|_{H^{\frac12}}.
\end{equation*}
\end{lem}
\begin{proof}
We use some symmetry property of the Kac's equation.
We write
the first term $I_{1a}=\frac12 I_{1a}+\frac12 I_{1a}$ and
we use the change of variables $\theta \to -\theta$.
We then have
\begin{equation}\label{I1a=}
I_{1a} = \int_{\RR_{\xi}}\int_{-\frac{\pi}{2}}^{\frac{\pi}{2}} \beta(\theta) \sin\theta \, \tilde A(\xi,\theta) \, 
  G_\delta(\xi \cos\theta) 
\hat g(\xi \cos\theta)
\overline {{\hat h}(\xi)} d\theta d\xi
\end{equation}
where
\begin{equation*}
\tilde A(\xi,\theta) = \frac12 \left( 
      \partial_\xi \hat f(  \xi \sin\theta ) 
    - \partial_\xi \hat f( -\xi \sin\theta )
  \right).
\end{equation*}
We compute
\begin{equation*}
\tilde A(\xi,\theta) =
\int_{\RR} v \sin(\xi v \sin\theta) f(v) dv 
\end{equation*}
and we estimate
\begin{equation*}
|\tilde A(\xi,\theta) | 
  \leq |\xi| \, |\sin\theta| \, \| f \|_{L^1_2}
  \leq \langle \xi \rangle \, |\sin\theta| \| f \|_{L^1_2}.
\end{equation*}
Finally we obtain
\begin{equation*}
 |I_{1a}| \leq C \| f \|_{L^1_2} 
  \|G_{\delta}g\|_{H^{\frac12}} \, \|h\|_{H^{\frac12}}.
\end{equation*}
\end{proof} 

\begin{lem}\label{lem I1a s=1/2}
Suppose that $s = \frac12$.
Then for any $0<{\alpha'}<1$, there exists a constant $C$ such that
\begin{equation*}
|I_{1a}| \leq
C \|f\|_{L^1_{1+{\alpha'}}}
\| G_\delta g \|_{H^{\frac{{\alpha'}}2}}
\| h \|_{H^{\frac{{\alpha'}}2}}.
\end{equation*}
\end{lem}
\begin{proof}
Following the proof of the previous lemma,
we consider again the identity \eqref{I1a=}
where
\begin{equation*}
\tilde A(\xi,\theta) =
\int_{\RR} v \sin(\xi v \sin\theta) f(v) dv.
\end{equation*}
We then estimate
\begin{equation*}
|\tilde A(\xi,\theta) | \leq |\xi|^{\alpha'} 
  |\sin\theta|^{\alpha'} \|f\|_{L^1_{1+{\alpha'}}}
  \leq \langle \xi \rangle^{\alpha'} 
   |\sin\theta|^{\alpha'} \|f\|_{L^1_{1+{\alpha'}}}.
\end{equation*}
Finally we obtain
\begin{equation*}
 |I_{1a}| \leq 
C \|f\|_{L^1_{1+{\alpha'}}}
\|G_{\delta}g\|_{H^{\frac{{\alpha'}}2}} \|h\|_{H^{\frac{{\alpha'}}2}}.
\end{equation*}
\end{proof} 

\begin{lem}\label{lem I1b}
There exists a constant $C$ such that
\begin{equation*}
|I_{1b}| \leq C 
\left(
  \|G_\delta f\|_{L^2_1} + \|G_\delta f\|_{H^{(\alpha-1)^+}}
\right)
\|G_\delta g\|_{H^{\frac{\alpha}2}_1}
\| h\|_{H^{\frac{\alpha}2}}.
\end{equation*}
\end{lem}
\begin{proof}
We estimate 
\begin{equation*}
I_{1b} = 
i\int_{\RR_{\xi}}\int_{-\frac{\pi}{2}}^{\frac{\pi}{2}}
\beta(\theta) \, \sin\theta \, 
(\partial_\xi \hat f)(\xi \sin\theta)
\left(G_\delta(\xi )- G_\delta(\xi\cos\theta ) 
\right) 
\hat g(\xi \cos\theta) \overline {{\hat h}(\xi)} \,
d\theta d\xi.
\end{equation*}
%1
By using lemma \ref{lem estim2 G},
\begin{align*}
|I_{1b}| 
\leq 
\int_{\RR_{\xi}}\int_{-\frac{\pi}{2}}^{\frac{\pi}{2}}
&\beta(\theta) \, \sin^2\textstyle{\frac{\theta}2} \, 
  |\sin\theta| \,
G_\delta(\xi\sin\theta)
\big| (\partial_\xi \hat f)(\xi\sin\theta) \big|
\\
&\langle \xi \rangle^{\frac{\alpha}2}
G_\delta(\xi\cos\theta ) |\hat g(\xi\cos\theta)|
\langle \xi \rangle^{\frac{\alpha}2} \overline {{\hat h}(\xi)} \,
d\theta d\xi.
\end{align*}
From
\begin{equation*}
 \|\langle \cdot \rangle^{\frac{\alpha}2} G_\delta \hat g\|_{L^\infty}
\leq
\|G_\delta g\|_{H^{\frac{\alpha}2}_1}
\end{equation*}
and the Cauchy-Schwarz inequality, we get
\begin{align*}
\int_{\RR_{\xi}}\int_{-\frac{\pi}{2}}^{\frac{\pi}{2}}
&\beta(\theta) \, \sin^2\textstyle{\frac{\theta}2} \, |\sin\theta| \,
G_\delta(\xi\sin\theta)
\big|(\partial_\xi \hat f)(\xi\sin\theta)\big|
\langle \xi \rangle^{\frac{\alpha}2} \overline {{\hat h}(\xi)} \,
d\theta d\xi
\\
\leq
&\left(
\int_{\RR_{\xi}}\int_{-\frac{\pi}{2}}^{\frac{\pi}{2}}
\beta(\theta) \, \sin^2\textstyle{\frac{\theta}2} \, |\sin\theta| \,
G_\delta(\xi\sin\theta)^2
\big|(\partial_\xi \hat f)(\xi\sin\theta)\big|^2 \,
d\theta d\xi
\right)^{1/2}
\\
&\left(
\int_{\RR_{\xi}}\int_{-\frac{\pi}{2}}^{\frac{\pi}{2}}
\beta(\theta) \, \sin^2\textstyle{\frac{\theta}2} \, |\sin\theta| \,
\langle \xi \rangle^{\frac{\alpha}2} \, |\hat h(\xi)|^2 \,
d\theta d\xi
\right)^{1/2}
\\
\leq
& \| G_\delta \partial_\xi \hat f \|_{L^2} \,\times\,
  \| \langle \cdot \rangle^{\frac{\alpha}2} \hat h \|_{L^2}.
\end{align*}
We then observe that from lemma \ref{lem estim1 G}
\begin{align*}
\| G_\delta (\partial_\xi \hat f) \|_{L^2}
&\leq
\| \partial_\xi (G_\delta \hat f) \|_{L^2}
+
\| (\partial_\xi G_\delta) \hat f \|_{L^2}
\\
&\leq
C\,
\| G_\delta f \|_{L^2_1} \,
\| \langle \cdot \rangle^{\alpha-1} \hat f \|_{L^2}
\end{align*}
and we conclude
\begin{equation*} 
|I_{1b}| \leq C
\left(
\| G_{\delta} f \|_{L^2_1} + \|G_{\delta} f\|_{H^{(\alpha-1)^+}}
\right) \,
\| G_{\delta} g \|_{H^{\frac{\alpha}2}_1} \,
\| h \|_{H^{\frac{\alpha}2}}.
\end{equation*}
\end{proof}

\medskip
{\bf Estimate of $I_2$.}
\medskip
We decompose $I_2 = I_{2a} + I_{2b}$
where
\begin{align*}
&I_{2a}
= i \int_{\RR_{\xi}}\int_{-\frac{\pi}{2}}^{\frac{\pi}{2}} \beta(\theta)
\hat f(\xi \sin\theta)
  \left(
  \cos\theta - 1
  \right) \,
G_\delta(\xi)
(\partial_\xi\hat g)(\xi\cos\theta)
\overline {{\hat h}(\xi)} \,d\theta d\xi,
\\
&I_{2b}
= i \int_{\RR_{\xi}}\int_{-\frac{\pi}{2}}^{\frac{\pi}{2}} \beta(\theta)
\hat f(\xi \sin\theta)
  \left(
  G_\delta(\xi) - G_\delta(\xi\cos\theta)
  \right) \, (\partial_\xi\hat g)(\xi\cos\theta)
\overline {{\hat h}(\xi)} \,d\theta d\xi.
\end{align*}
\begin{lem}\label{lem I2a}
There exists a constant $C$ such that
\begin{equation*}
|I_{2a}| \leq
C \, \| G_{\delta} f \|_{L^2_1}
\| G_{\delta} g \|_{H^{(\alpha-1)^+}}
\|h\|_{L^2}.
\end{equation*}
\end{lem}

\begin{proof}
For $I_{2a}$ we use \eqref{G} and 
$\cos\theta -1 = -2\sin^2\textstyle{\frac{\theta}2}$ :
\begin{align*}
I_{2a}
&\leq 
C \|G_{\delta}f\|_{L^2_1}
\|G_{\delta}(vg)\|_{L^2}
\|h\|_{L^2}
\\
&\leq
C \,
\| G_{\delta} f \|_{L^2_1}
\left(
\| G_{\delta} g \|_{L^2}
+
\| G_{\delta} g \|_{H^{(\alpha-1)^+}}
\right)
\| h \|_{L^2}.
\end{align*}
\end{proof} 

\begin{lem}\label{lem I2b}
There exists a constant $C$ such that
\begin{equation*}
|I_{2b}| \leq 
C \, 
\| G_{\delta}f \|_{L^2_1}
\left(
\| G_{\delta} g \|_{H^{\frac{\alpha}2}_1}
+
\| G_{\delta} g \|_{H^{(\frac{3\alpha}2-1)^+}}
\right)
\| h \|_{H^{\frac{\alpha}2}}.
\end{equation*}
\end{lem}
\begin{proof}
Using lemma \ref{lem estim2 G} we get
\begin{equation*}
I_{2b}
\leq 
C \|G_{\delta} f\|_{L^\infty} \,
\|\langle \cdot \rangle^{\frac{\alpha}2} 
  G_{\delta} (\partial_\xi \hat g) \|_{L^2}
\|\langle \cdot \rangle^{\frac{\alpha}2} \hat h\|_{L^2}
\end{equation*}
and 
\begin{equation*}
\| \langle \cdot \rangle^{\frac{\alpha}2} G_{\delta} (\partial_\xi \hat g) \|_{L^2}
\leq
\| \langle \cdot \rangle^{\frac{\alpha}2} \partial_\xi (G_{\delta} \hat g) \|_{L^2}
+
\|\langle \cdot \rangle^{\frac{\alpha}2} (\partial_\xi G_{\delta}) \hat g \|_{L^2}.
\end{equation*}
\end{proof} 

\medskip
{\bf Estimate of $I_3$ }
\medskip
We recall
\begin{equation*}
I_3 = 
i \int_{\RR_{\xi}}\int_{-\frac{\pi}{2}}^{\frac{\pi}{2}} \beta(\theta) 
\hat f(\xi \sin\theta)
  \left\{
  \partial_\xi G_\delta (\xi) 
 -(\partial_\xi G_\delta) (\xi\cos\theta)
  \right\}
\, \hat g(\xi\cos\theta) \,
\overline {{\hat h}(\xi)} \,d\theta d\xi.
\end{equation*}
\begin{lem}\label{lem I3}
There exists a constant $C$ such that
\begin{equation*}
I_3 
\leq
C \, 
\| G_{\delta}f \|_{L^2_1}
\| G_{\delta} g \|_{H^{(\alpha-\frac12)}}
\| h \|_{H^{(\alpha-\frac12)}}.
\end{equation*}
\end{lem}
\begin{proof}
\begin{equation*}
I_3 
\leq
C \, \| G_{\delta} \hat f \|_{L^\infty}
\| \langle \cdot \rangle^{(\alpha-\frac12)} G_{\delta} \hat g \|_{L^2}
\| \langle \cdot \rangle^{(\alpha-\frac12)} \hat h \|_{L^2}.
\end{equation*}
\end{proof} 

\noindent
{\bf Proof of Proposition \ref{Prop estim com1}.}
We use the previous lemmas \ref{lem I1a} and \ref{lem I1b}-\ref{lem I3}.
By summing the above estimates, we deduce from \eqref{com1<=}
\begin{align*}
\big| \big( (v G_\delta) &K(f, g) , h\big)_{L^2}
  - \big(K(f, (v G_\delta) g) , h\big)_{L^2} \big|
\\
&\leq\,
C \|f\|_{L^1_2}
\|G_\delta g \|_{H^{\frac12}} \|h\|_{H^{\frac12}}
+ C \|G_\delta f\|_{H^{(\alpha-1)^+}_1}
\|G_\delta g \|_{H^{\frac{\alpha}2}}
\|h\|_{H^{\frac{\alpha}2}}
\\
&+ C \|G_\delta f\|_{L^2_1}
\|G_\delta g \|_{H^{\left(\frac{3\alpha}2-1\right)^+}_1}
\|h\|_{H^{\frac{\alpha}2}}
+ C \|G_\delta f\|_{L^2_1}
\|G_\delta g \|_{H^{\left(\alpha-\frac12\right)^+}}
\|h\|_{H^{\left(\alpha-\frac12\right)^+}}.
\end{align*}
Taking $\alpha=1$,
this finishes the proof of Proposition \ref{Prop estim com1}.
\hfill$\square$

\noindent
{\bf Proof of Proposition \ref{Prop estim com1 s=1/2}.}
We recall $s=\frac12$.
We have from \eqref{com1<=}
\begin{equation*}
\big|\big( (v G_\delta) K(f, g) , h\big)_{L^2} 
  - \big(K(f, (v G_\delta) g) , h\big)_{L^2} \big|
\leq |I_{1a}| + |I_{1b}| + 
|I_{2a}| + |I_{2b}| + |I_3|.
\end{equation*}
We use the lemma \ref{lem I1a s=1/2}
and the lemmas \ref{lem I1b}-\ref{lem I3}
taking $0<\alpha<1$,
and this concludes the proof.
\hfill$\square$

We now estimate some scalar product terms
which involve the derivative 
of the mollifier with respect to time:
\begin{lem} \label{lem temp}
There exists $C>0$ such that
\begin{equation}\label{estim time0}
\left|\left(
  (\partial_{t} G_{\delta})(t,D_v)f(t,\cdot) \,,\,
  G_{\delta}(t,D_v)f(t,\cdot)
\right)_{L^2} \right|
\leq
  C \, \|G_{\delta}f\|^2_{H^{\alpha/2}},
\end{equation}
and
\begin{equation}\label{estim time1}
\left|\left(
  v(\partial_{t} G_{\delta})(t,D_v)f(t,\cdot) \,,\,
  vG_{\delta}(t,D_v)f(t,\cdot)
\right)_{L^2} \right|
\leq
  C \,
\Big(
\|G_{\delta}f \|^2_{H^{\alpha/2}_1}
+
\|G_{\delta}f \|^2_{H^{\alpha-\frac12}_1}
\Big).
\end{equation}
\end{lem}
\begin{proof}
We have by the Plancherel formula
\begin{equation*}
\left((\partial_{t} G_{\delta})(t,D_v)f(t,\cdot) \,,\,
G_{\delta}(t,D_v)f(t,\cdot)\right)_{L^2}
=
\int
(\partial_t G_\delta) \hat f \,
\overline{G_\delta \hat f} d\xi.
\end{equation*}
The estimate \eqref{estim time0} can be deduced directly from 
lemma \ref{lem estim1 G}.
For \eqref{estim time1}, we compute
\begin{align*}
\big(
  v(\partial_{t} & G_{\delta})(t,D_v)f(t,\cdot) \,,\,
  vG_{\delta}(t,D_v)f(t,\cdot)
\big)_{L^2}
\\
&=
\int
\left\{
\partial_\xi
\left(
  \frac{c_0 \langle \xi \rangle^{\alpha}}{1+\delta e^{c_0t\langle \xi \rangle^{\alpha}}}
\right)
 (G_{\delta} \hat f) \,
+
\left(
  \frac{c_0 \langle \xi \rangle^{\alpha}}{1+\delta e^{c_0 t \langle \xi \rangle^{\alpha}}}
\right)
\partial_\xi (G_{\delta} \hat f) \,
\right\}
\,
\overline{\partial_\xi (G_\delta \hat f) }
d\xi
\end{align*}
and we use the following estimate
\begin{equation*}
\left|\partial_{\xi}\left(
   \frac{\langle \xi \rangle^{\alpha}}{1+\delta e^{c_0 t \langle \xi \rangle^{\alpha}}}
\right) \right|
\leq
C\langle \xi \rangle^{2\alpha-1}.
\end{equation*}
\end{proof} 

\renewcommand{\theequation}{\thesection.\arabic{equation}}
\setcounter{equation}{0}
\section{SOBOLEV REGULARIZING EFFECT FOR KAC'S EQUATION}\label{section Sobolev}
In this section, we prove the regularity in weighted Sobolev spaces 
of the weak solutions for the Cauchy problem of the Kac's equation.
\begin{theo}\label{Th Sobolev}
Assume that the initial datum 
$f_0\in L^1_{2+2s}\bigcap L\log L(\RR)$,
and the cross-section weak $\beta$ satisfies \eqref{beta}
with $\frac12 \leq s<1$.
If $f\in L^{\infty}( ]0,+\infty[ ; L^1_{2+2s} \bigcap L\log L(\RR) )$
is a nonnegative weak solution 
of the Cauchy problem \eqref{eq Kac},
then $f(t,\cdot)\in H^{+\infty}_2(\RR)$ for any $t>0$.
\end{theo}
\noindent
{\bf Remark.}
This Theorem has been proved in \cite{LX09}
in the case $0<s<\frac12$.

We also obtain the following propagation of Sobolev regularity:

\begin{corol}\label{Corol Sobolev}
Under the assumptions  of Theorem \ref{Th Sobolev},
for any $T_0>0$, 
there exists a constant $C$
which depends only on $\beta$ and
$\| f \|_{  L^{\infty}( ]0,+\infty[ ; L^1_{2+2s} \bigcap L\log L(\RR) ) }$
such that
\begin{equation*}
\forall t \geq T_0, \quad
\| f(t,\cdot) \|_{H^2_2}
\leq
 e^{C (t-T_0)} 
\| f(T_0,\cdot) \|_{H^2_2}.
\end{equation*}
\end{corol}

Throughout this section, 
we will distinguish the case $\frac12<s<1$ and the limit case $s=\frac12$.
We introduce as in \cite{MUXY09} the mollifier of polynomial type
\begin{equation*}
M_{\delta}(t,\xi)= \frac{\langle \xi \rangle^{Nt-1}}{(1+\delta|\xi|^2)^{N_0}}
\end{equation*}
for $0<\delta<1$, $t\in [0,T_0]$ and $2N_0= T_0N+4$.
\begin{lem}\label{lem estime M}
We have that for any $0<\delta<1$ and $0\leq t\leq T_0$, $\xi\in \RR$,
\begin{align*}
\left| \partial_{t} M_{\delta}(t,\xi)\right| 
  \leq N\log(\langle \xi \rangle) M_{\delta}(t,\xi).
\end{align*}
For $-\frac{\pi}{4}\leq \theta \leq \frac{\pi}{4}$,
\begin{align*}
\left| M_\delta(t, \xi) - M_\delta(t, \xi \cos \theta) \right|
  & \leq
    C \sin^2 \textstyle{\frac{\theta}2}M_\delta(t, \xi \cos \theta),
\\
  \left| (\partial_{\xi} M_\delta)(t, \xi)
    - (\partial_{\xi} M_\delta)(t, \xi \cos \theta) \right|
  &\leq
    C \sin^2 \textstyle{\frac{\theta}2}
    \langle \xi \rangle^{-1 } M_\delta(t, \xi \cos \theta),
\\
  \left|(\partial^2_{\xi} M_\delta)(t, \xi)
  - (\partial^2_{\xi} M_\delta)(t, \xi \cos \theta) \right|
  &\leq
    C \sin^2 \textstyle{\frac{\theta}2}\langle \xi
    \rangle^{-2 } M_\delta(t, \xi \cos \theta).
\end{align*}
\end{lem}
\begin{proof}
We compute
\begin{align*}
\log M_\delta(t, \xi) 
&= 
\frac{N t -1}{2} \log (1+\xi^2)
- 
N_0 \log (1+\delta \xi^2),
\\
\partial_t M_\delta(t, \xi)
&= 
\frac{N}{2} \log (1+\xi^2) \, M_\delta(t, \xi).
\end{align*}
Using the estimates
\begin{align*}
\left| \partial^{k}_{\xi} (M_\delta(t,\xi))\right|
&\leq 
    C_k\langle \xi \rangle^{-k} M_\delta(t,\xi),
\\
\left| M_\delta(t,\xi)\right|
&\leq 
    C M_\delta(t,\xi \cos \theta)
\end{align*}
and the Taylor formula, we obtain the proof of the lemma.
\end{proof} 
We estimate the first commutator:
\begin{prop}\label{Prop estim comM0}
Let $f,g \in L^2_1$ and $h\in L^2(\RR)$,
then we have that 
\begin{equation*}
\big| \big( M_\delta K(f, g) , h \big)_{L^2}
  - \big( K(f, M_\delta g) , h \big)_{L^2} \big|
\leq C \, \| f \|_{L^1} \|M_{\delta} g\|_{L^2} \|h\|_{L^2}.
\end{equation*}
\end{prop}
\begin{proof}
By the definition of $M_\delta$, we have
\begin{align*}
\cF(M_{\delta}f)(\xi) &= M_{\delta}\hat{f}(\xi),
\\
\cF(v M_{\delta}f)(\xi) &=i\partial _{\xi}
  \left( M_{\delta}(t,\xi)\hat{f}(\xi)\right).
\end{align*}
We now use the Bobylev formula \eqref{bobylev}
and the Plancherel formula
\begin{align*}
\big( &M_\delta K(f, g) , h \big)_{L^2} 
- \big( K(f, M_\delta g) , h \big)_{L^2}
\\
&= \int_{\RR_{\xi}} \int_{-\frac{\pi}{2}}^{\frac{\pi}{2}}
  \beta(\theta)\hat{f}(\xi \sin \theta)
  \left\lbrace M_{\delta}(\xi) 
  - M_{\delta}(\xi \cos \theta) \right\rbrace
  \hat{g}(\xi\cos\theta)\overline{\hat{h}(\xi)} d\theta d\xi.
\end{align*}
By the previous formula, lemma \ref{lem estime M},
the Cauchy-Schwarz inequality and \eqref{cond beta} we have 
\begin{align*}
\big|\big( M_\delta K(f, g) , h\big)_{L^2} 
  &- \big(K(f, M_\delta g) , h\big)_{L^2} \big|
\\
&\leq 
\int_{\RR_{\xi}}\int_{-\frac{\pi}{2}}^{\frac{\pi}{2}}
\beta(\theta)\sin^2\textstyle{\frac{\theta}2}
|\hat{f}(\xi\sin\theta)|
\times
|M_{\delta}(\xi\cos\theta)\hat{g}(\xi\cos\theta) \hat{h}(\xi)|
d\theta d\xi
\\
&\leq
C \|\hat f\|_{L^{\infty}}
  \|M_{\delta} g\|_{L^2} \|h\|_{L^2}
\\
&\leq
C \, \| f \|_{L^1} \|M_{\delta} g\|_{L^2} \|h\|_{L^2}.
\end{align*}
\end{proof} 
In the same spirit of Proposition \ref{Prop estim com1},
we will use some symmetry property of the Kac's equation
to estimate the weighted commutator.

\begin{prop}\label{Prop estim comM2}
Suppose that $\frac12< s < 1$. We then have :
\begin{equation*}
\left|
\big(  (v^2 M_\delta) K(f, g) , h\big)_{L^2} 
- \big(K(f, (v^2 M_\delta) g) , h\big)_{L^2}
\right|
\leq
  C \, \|f\|_{L^1_2} \,
  \|M_\delta g\|_{H^{\frac12}_2} \, 
  \|h\|_{H^{\frac12}}.
\end{equation*}
\end{prop}
\begin{proof}
We have
\begin{align*}
&\big((v^2 M_\delta) K(f, g) , h\big)_{L^2} 
- \big(K(f, (v^2 M_\delta) g) , h\big)_{L^2} 
\\
&=
-\int\int_{-\frac{\pi}{2}}^{\frac{\pi}{2}}
\beta(\theta) (\sin^2\theta)
(\partial^2_{\xi}\hat{f})(\xi \sin \theta)
M_{\delta}(\xi)
\hat{g}(\xi\cos\theta)\overline{\hat{h}(\xi)}
d\theta d\xi
\\
&- 2 \int\int_{-\frac{\pi}{2}}^{\frac{\pi}2}
\beta(\theta)
\sin\theta (\partial_{\xi}\hat{f})(\xi \sin \theta)
\,
\partial_{\xi}
\Bigl(
  M_{\delta}(\xi) \, \hat{g}(\xi\cos\theta)
\Bigr) \,
\overline{\hat{h}(\xi)}
d\theta d\xi
\\
&- \int\int_{-\frac{\pi}{2}}^{\frac{\pi}{2}}
\beta(\theta)\hat{f}(\xi \sin \theta)
\left\lbrace
  \partial^2_{\xi}
  \Bigl(
    M_\delta(t, \xi) \, \hat g(\xi \cos \theta)
  \Bigr)
  -
  \partial^2_{\xi}
  \Bigl(
      M_\delta \, \hat g
  \Bigr)
     (\xi \cos \theta)
\right\rbrace
\overline{\hat{h}(\xi)}d\theta d\xi
\\
&= B_1 + B_2 + B_3.
\end{align*}
Then
\begin{equation*}
|B_1| 
\leq C \, \| f \|_{L^1_2} \, \|M_{\delta}g\|_{L^2} \, \|h\|_{L^2}.
\end{equation*}
For $B_2$, we will use the symmetry and
the change of variables $\theta \to - \theta$
(see proof of lemma \ref{lem I1a}).
We write $B_2 = B_{2a} + B_{2b}$ where
\begin{align*}
B_{2a} &=
- 2 \int_{\RR_{\xi}}\int_{-\frac{\pi}{2}}^{\frac{\pi}2}
\beta(\theta)
\sin\theta (\partial_{\xi}\hat{f})(\xi \sin \theta)
\,
\left( \partial_{\xi} M_{\delta} \right)(\xi) \,
\hat{g}(\xi\cos\theta)
\overline{\hat{h}(\xi)}
d\theta d\xi,
\\
B_{2b} &=
- 2 \int_{\RR_{\xi}}\int_{-\frac{\pi}{2}}^{\frac{\pi}2}
\beta(\theta)
\sin\theta (\partial_{\xi}\hat{f})(\xi \sin \theta)
\,
\left( \partial_{\xi} \hat{g} \right) (\xi\cos\theta)\,
\cos\theta \,
M_{\delta}(\xi)
\overline{\hat{h}(\xi)}
d\theta d\xi.
\end{align*}
The symmetry and the estimate of lemma \ref{lem estime M}
implies
\begin{align*}
|B_{2a}| &\leq 
C \int_{\RR_{\xi}}\int_{-\frac{\pi}{2}}^{\frac{\pi}2}
\beta(\theta) 
(\sin^2\theta)  |\xi|
\| \partial_{\xi}^2 \hat f \|_{L^\infty}
\langle \xi \rangle^{-1} M_{\delta}(\xi\cos\theta) 
 \,
|\hat{g}(\xi\cos\theta)|
\times
|\hat h(\xi)|
d\theta d\xi
\\
&\leq C \| f \|_{L^1_2} 
\|M_\delta g\|_{L^2} \|h\|_{L^2}.
\end{align*}
We note that
\begin{equation*}
\|\langle \xi \rangle^{\frac12} \, M_\delta \partial_\xi \hat g \|_{L^2}
=
\| \, M_\delta (v g) \|_{H^{\frac12}}
\leq C \|M_\delta g \|_{H^{\frac12}_1}.
\end{equation*}
Using again the symmetry
and the previous estimate,
we get
\begin{align*}
|B_{2b}| &\leq
C \int_{\RR_{\xi}}\int_{-\frac{\pi}{2}}^{\frac{\pi}2}
\beta(\theta) 
\sin^2\theta  \langle \xi \rangle
\|\partial_{\xi}^2 \hat{f} \|_{L^\infty}
\,
| M_{\delta}(\xi\cos\theta) 
\partial_{\xi}
  \hat{g}(\xi\cos\theta)|
\times
|\hat{h}(\xi)|
d\theta d\xi
\\
&\leq
C\, \| f \|_{L^1_2} \, \|M_\delta g\|_{H^{\frac12}_1}
    \|h\|_{H^{\frac12}}.
\end{align*}
For $B_3$ we have 
\begin{align*}
B_3 =
- \int_{\RR_{\xi}}\int_{-\frac{\pi}{2}}^{\frac{\pi}{2}}
&\beta(\theta) 
\hat{f}(\xi \sin \theta)
\\
&\left\lbrace
  \partial^2_{\xi}
  \Bigl(
    M_\delta(t, \xi) \, \hat g(\xi \cos \theta) 
  \Bigr)
  -
  \partial^2_{\xi}
  \Bigl(
      M_\delta \, \hat g
  \Bigr)
     (\xi \cos \theta) 
\right\rbrace 
\overline{\hat{h}(\xi)}d\theta d\xi
\end{align*}
and we compute
\begin{equation*}
\partial^2_{\xi}
\Bigl(
  M_\delta(t, \xi) \, \hat g(\xi\cos\theta) 
\Bigr)
-
\partial^2_{\xi}
\Bigl(
  M_\delta \, \hat g
\Bigr)(\xi \cos \theta) 
=
D_1 + D_2 + D_{3} + D_{4} + D_{5}
\end{equation*}
where
\begin{align*}
D_1 &= 
\left(
    \partial_\xi^2 M_\delta(\xi) 
  - \partial_\xi^2 M_\delta(\xi\cos\theta)
\right) \,
\hat g(\xi\cos\theta),
\\
D_2 &=
\left(
     M_\delta(\xi) 
  -  M_\delta(\xi\cos\theta)
\right) \,
\partial_\xi^2\hat g(\xi\cos\theta),
\\
D_3 &=
(\cos^2\theta -1) M_\delta(\xi\cos\theta)
\partial_\xi^2\hat g(\xi\cos\theta),
\\
D_4 &=
2\left(
     \partial_\xi M_\delta(\xi) 
  -  \partial_\xi M_\delta(\xi\cos\theta)
\right) \,
\partial_\xi \hat g(\xi\cos\theta),
\\
D_5 &=
2(\cos\theta - 1)
\partial_\xi M_\delta(\xi) \,
\partial_\xi \hat g(\xi\cos\theta).
\end{align*}
For $1\leq i \leq 5$, we note
$\displaystyle J_i = 
- \int_{\RR_{\xi}}\int_{-\frac{\pi}{2}}^{\frac{\pi}{2}}
\beta(\theta)\hat{f})(\xi \sin \theta)
D_i(\xi,\theta)
\overline{\hat{h}(\xi)}d\theta d\xi
$.
\\
We successively estimate :
\begin{align*}
|J_1| 
&\leq 
C \, \| \hat f \|_{L^\infty}
\int_{\RR_{\xi}} \int_{-\frac{\pi}{2}}^{\frac{\pi}{2}}
\beta(\theta) \textstyle{\sin^2\frac\theta2}
   M_\delta(\xi\cos\theta) |\hat g(\xi\cos\theta)\,
\hat h(\xi)| d\theta d\xi
\\
&\leq
C \, \|f\|_{L^1} \, \|M_\delta g\|_{L^2} \|h\|_{L^2} ,
\\
|J_2| + |J_3|
&\leq
C \, \|\hat f \|_{L^\infty}
\int_{\RR_{\xi}} \int_{-\frac{\pi}{2}}^{\frac{\pi}{2}}
\beta(\theta) 
\textstyle{(\sin^2\frac\theta2 + \sin^2\theta)}
   M_\delta(\xi\cos\theta) |\partial_\xi^2\hat g(\xi\cos\theta) \,
\hat h(\xi)| d\theta d\xi
\\
&\leq
C \, \|f\|_{L^1} \, \|M_\delta g\|_{L^2_2} \|h\|_{L^2},
\\
|J_4| + |J_5| 
&\leq
C \, \| \hat f \|_{L^\infty}
\int_{\RR_{\xi}} \int_{-\frac{\pi}{2}}^{\frac{\pi}{2}}
\beta(\theta) 
\textstyle{\sin^2\frac\theta2)}
   M_\delta(\xi\cos\theta) |\partial_\xi \hat g(\xi\cos\theta) \,
 \hat h(\xi)| d\theta d\xi
\\
&\leq
C \, \|f\|_{L^1} \, \|M_\delta g\|_{L^2_1} \|h\|_{L^2}.
\end{align*}
From the previous inequalities we deduce
\begin{equation*}
|B_3| \leq 
C\,\|f\|_{L^1} \, \|M_{\delta }g\|_{L^2_2} \, \|h\|_{L^2} 
\end{equation*}
and this finishes the proof
of the Proposition \ref{Prop estim comM2}.
\end{proof}

For the case $s=\frac12$, we will need a different
estimate of the weighted commutator.
\begin{prop}\label{Prop estim comM2 s=1/2}
Assume that $s =\frac12$.
Then for any $0<{\alpha'}<1$ we have
\begin{equation}\label{com2M s=1/2}
\left|
\big( (v^2 M_\delta) K(f, g) , h\big)_{L^2}
- \big(K(f, (v^2 M_\delta) g) , h\big)_{L^2}
\right|
\leq
  C \, \|f\|_{L^1_2} \,
  \|M_\delta g\|_{H^{\frac{{\alpha'}}2}_2} \,
  \|h\|_{H^{\frac{{\alpha'}}2}}.
\end{equation}
\end{prop}

The proof of this Proposition use the same arguments of
Proposition \ref{Prop estim comM2}
and lemma \ref{lem I1a s=1/2}.

{\bf Proof of the Theorem \ref{Th Sobolev}}.

\noindent
{\bf - Case :} $\frac12<s<1$.

We consider $f\in L^1_{2+2s}\cap L\log L$
a weak solution of the Cauchy problem \eqref{eq Kac}
and we multiply the equation with the test function
\begin{equation*}
\varphi(t,v) = M_\delta(t,D_v) (1+v^4) M_\delta(t,D_v) f(t,v).
\end{equation*}
Therefore we obtain the equality
\begin{equation}\label{dt f phi}
(\partial_t f, \varphi)_{L^2}
=
(K(f,f), \varphi)_{L^2}.
\end{equation}
Using some similar arguments in \cite{MUXY09},
we can suppose that
$\varphi\in C^1([0,T_0];H^{5}_{-2+2s}(\RR))$.
We compute
\begin{align*}
(M_\delta\partial_t f, M_\delta f)_{L^2}
&+
(v^2 M_\delta\partial_t f, v^2 M_\delta f)_{L^2}
\\
&=
(M_\delta K(f,f), M_\delta f)_{L^2}
+
(v^2 M_\delta K(f,f), v^2 M_\delta f)_{L^2}.
\end{align*}
We will use the following notations
\begin{align*}
\hbox{time}_0
&= \left((\partial_t M_\delta) f,
   M_\delta f\right)_{L^2},
\\
\hbox{time}_2
&= \left(v^2 (\partial_t M_\delta) f,
   v^2 M_\delta f\right)_{L^2},
\end{align*}
and, concerning the commutators of the Kac operator
and the weighted mollifier,
\begin{align*}
\hbox{com}_0
&= \left( M_\delta K(f,f), M_\delta f\right)_{L^2}
 - \left( K(f,M_\delta f), M_\delta f\right)_{L^2},
\\
\hbox{com}_2
&= \left( v^2 M_\delta K(f,f), v^2 M_\delta f\right)_{L^2}
 - \left( K(f,v^2 M_\delta f), v^2 M_\delta f\right)_{L^2}.
\end{align*}
Therefore the relation \eqref{dt f phi} become
\begin{align*}
\frac12 \frac{d}{dt}
\Bigl( &\|M f\|_{L^2}^2 + \|v^2 M f\|_{L^2}^2 \Bigr)
- (K(f,M_\delta f) \,,\, M_\delta f)_{L^2}
- (K(f , v^2 M_\delta f) \,,\, v^2 M_\delta f)_{L^2}
\\
&=
\hbox{time}_0 + \hbox{time}_2
+\hbox{com}_0 + \hbox{com}_2.
\end{align*}
From the coercivity inequality of Proposition \ref{Prop coer},
we derive the following differential inequation
\begin{align}\label{ineq energ}
\frac12 \frac{d}{dt}
\Bigl( \| M f \|_{L^2}^2 + \| v^2 M f \|_{L^2}^2 \Bigr)
  &+
c_f \| M f \|_{H^s_2}^2
\\
&\leq
\hbox{time}_0 + \hbox{time}_2 +\hbox{com}_0 + \hbox{com}_2
+ C \| f \|_{L^1} \| M_\delta f \|_{L^2_2}. \nonumber
\end{align}
\begin{lem}\label{lem estim temp}
Assume that $0< s < 1$ and $\varepsilon>0$.
Then there exists a constant $C_\varepsilon$
such that :
\begin{align*}
\left|
  \left((\partial_t M_\delta) f,
   M_\delta f\right)_{L^2}
\right|
&\leq
\varepsilon \, \|M_\delta f\|^2_{H^{s}_2} \,
   + C_\varepsilon \|M_\delta f\|^2_{L^2},
\\
\left|
  \left(v^2 (\partial_t M_\delta) f,
   v^2 M_\delta f\right)_{L^2}
\right|
&\leq
  \varepsilon \, \|M_\delta f\|^2_{H^{s}_2} \,
   + C_\varepsilon \|M_\delta f\|^2_{L^2}.
\end{align*}
\end{lem}
\begin{proof}
We compute
\begin{equation*}
\hbox{time}_0
=
\left(
  \frac{N}2 \, \log(1+\xi^2) \, M_\delta \hat f \,,\,
  M_\delta \hat f
\right)_{L^2}.
\end{equation*}
For $\varepsilon>0$, there exists a constant $C_\varepsilon$
such that :
\begin{equation}\label{log eps}
\frac{N}2 \log (1+\xi^2) \, 
\leq \,
\varepsilon (1+\xi^2)^{s} + C_\varepsilon.
\end{equation}
Therefore 
\begin{equation*}
|\hbox{time}_0|
\leq 
  \varepsilon
  \left(
   (1+\xi^2)^{s} M_\delta \hat f \,,\,
    M_\delta \hat f 
  \right)_{L^2}
+
  C_\varepsilon
  \left(
    M_\delta \hat f \,,\,
    M_\delta \hat f 
  \right)_{L^2}.
\end{equation*}
We estimate the term 
$\hbox{time}_2
=
\left(
  \partial_\xi^2
  \left(
    \partial_t M_\delta \hat f
  \right)
 \, , \,
\cF (v^2 M_\delta f)
\right)_{L^2}
$.
We compute
\begin{align*}
\partial_\xi^2
  \left(
    \partial_t M_\delta \hat f 
  \right)
&=
\underbrace{ 
\partial_\xi^2
  \left(
    \frac{N}2 \log(1+\xi^2)
  \right)
}_{\leq C}
\,
M_\delta \hat f
+
\underbrace{
2 
\partial_\xi
  \left(
    \frac{N}2 \log(1+\xi^2)
  \right)
}_{\leq C}
\,
\partial_\xi \left( M_\delta \hat f \right)
\\
&+
\frac{N}2 \log (1+\xi^2) \,
\partial_\xi^2 \left( M_\delta \hat f \right).
\end{align*}
Using again the inequality \eqref{log eps},
\begin{align*}
|\hbox{time}_2|
&\leq 
C
\left|\left(
  M_\delta \hat f \, 
  ,\, 
  \cF (v^2 M_\delta f)
\right)_{L^2}\right|
+
C
\left|\left(
  \partial_\xi
  \left(M_\delta \hat f\right) \, 
  ,\, 
  \cF (v^2 M_\delta f)
\right)_{L^2}\right|
\\
&+
\varepsilon
\left|\left(
  \partial_\xi^2
  \left(M_\delta \hat f\right) \, 
  ,\, 
  \cF (v^2 M_\delta f)
\right)_{L^2}\right|
+
C_\varepsilon
\left|\left(
  (1+\xi^2)^{s}
  \partial_\xi^2
  \left(M_\delta \hat f\right) \, 
  ,\, 
  \cF (v^2 M_\delta f)
\right)_{L^2}\right|
\end{align*}
where
\begin{align*}
\left|\left(
  M_\delta \hat f \, 
  ,\, 
  \cF {(v^2 M_\delta f)}
\right)_{L^2}\right|
&=
\| v M_\delta f \|_{L^2_2}^2
\leq
\| M_\delta f \|_{L^2_2}^2,
\\
\left|\left(
  \partial_\xi (M_\delta \hat f) \, 
  ,\, 
  \cF (v^2 M_\delta f)
\right)_{L^2}\right|
&= 
\left|\left(
  v M_\delta f \, 
  ,\, 
   v^2 M_\delta f
\right)_{L^2}\right|
\leq
\| M_\delta f \|_{L^2_2},
\\
\left|\left(
  \partial_\xi^2 (M_\delta \hat f) \, 
  ,\, 
  \cF (v^2 M_\delta f)
\right)_{L^2}\right|
&= 
\left|\left(
  v^2 M_\delta f \, 
  ,\, 
   v^2 M_\delta f
\right)_{L^2}\right|
\leq
\| M_\delta f \|_{L^2_2},
\\
\left|\left(
  (1+\xi^2)^{s}
  \partial_\xi^2
  \left(M_\delta \hat f\right) \, 
  ,\, 
  \cF (v^2 M_\delta f)
\right)_{L^2}\right|
&=
\| M_\delta f \|_{H^{s}_2}^2.
\end{align*}
This concludes the proof of lemma \ref{lem estim temp}.
\end{proof} 
Plugging the estimates of Propositions 
\ref{Prop estim comM0}, \ref{Prop estim comM2}
and lemma \ref{lem estim temp}
into \eqref{ineq energ}, we get
\begin{align*}
\frac12 \frac{d}{dt}
&\Big( 
  \|M_\delta f\|_{L^2}^2 + \|v^2 M_\delta f\|_{L^2}^2 
 \Big)
  + c_f \|M_\delta f\|_{H^s_2}^2 
\\
&\leq
 \varepsilon \|M_\delta f\|_{H^{s}_2}^2 
+ C_\varepsilon \|f\|_{L^2_2}^2 
+ 
  C \|f\|_{L^1_2} \|M_\delta f\|_{H^{\frac12}_2}^2.
\end{align*}
From the interpolation estimate
\begin{equation}\label{interp H1/2}
\| g \|_{H^{\frac12}}^2 
\leq
\lambda \| g \|_{H^s}^2
+
\lambda^{\frac{-1}{2s-1}} \| g \|_{L^2}^2,
\end{equation}
we deduce
\begin{align*}
\frac12 \frac{d}{dt}
&\Big( \| M_\delta f \|_{L^2}^2 + \|v^2 M_\delta f\|_{L^2}^2 \Big)
  +
c_f \|M_\delta f\|_{H^s_2}^2 
\\
&\leq
\left(
 \varepsilon + C \lambda \| f \|_{L^1_2}
\right)
  \,  \|M_\delta f\|_{H^{s}_2}^2
+
\left(
    C_\varepsilon 
  + C \lambda^{\frac{-1}{2s-1}} \| f \|_{L^1_2}
\right) \,
\|M_\delta f\|_{L^2_2}^2.
\end{align*}
Choosing $\varepsilon$ and $\lambda$ small enough,
we get
\begin{equation}\label{dt Mf <}
\frac12 \frac{d}{dt}
\Big( \|M_\delta f\|_{L^2}^2 + \|v^2 M_\delta f\|_{L^2}^2 \Big)
\leq
C \|M_\delta f\|_{L^2_2}^2 
\leq
C\Big( \|M_\delta f\|_{L^2}^2 + \|v^2 M_\delta f\|_{L^2}^2 \Big).
\end{equation}
From Gronwall's lemma we have
\begin{equation*}
\|M_\delta f\|_{L^2}^2 + \|v^2 M_\delta f\|_{L^2}^2 
\leq
e^{2 C t}
\left(
\|M_\delta(0) f_0\|_{L^2}^2 + \|v^2 M_\delta(0) f_0\|_{L^2}^2
\right),
\end{equation*}
that is
\begin{equation*}
\|M_\delta f\|_{L^2_2}^2
\leq
 e^{2 C t} 
 \|M_\delta(0) f_0\|_{L^2_2}^2.
\end{equation*}
We write :
\begin{equation*}
\| (1-\delta \Delta)^{-N_0} f \|_{H^{N t -1}_2}^2
\leq
 e^{2 C t} 
\| (1-\delta \Delta)^{-N_0} f_0 \|_{H^{-1}_2}^2.
\end{equation*}
By Fatou's lemmas, letting $\delta\to0$,%
\begin{equation*}
\| f \|_{H^{N t -1}_2}^2 
\leq
 e^{2 C t} 
\| f_0 \|_{H^{-1}_2}^2 
\leq
C' e^{2 C t} 
\| f_0 \|_{L^1_2}^2.
\end{equation*}
For $t\in[0, T_0]$ we have proved
\begin{equation*}
  (1 + |D_v|^2)^m f(t, \cdot) \in L_2^2(\RR)
\end{equation*}
for all $T_0>0$ and $m=N t - 1>0$.
Therefore we have obtained that $f(t,\cdot)\in H_2^m(\RR)$
and that concludes
the proof of Theorem \ref{Th Sobolev} in the case $\frac12<s<1$.

\noindent
{\bf - Case} $s=\frac12$.
The proof is similar to the case $\frac12<s<1$.
We choose $\alpha'=\frac12$ in 
Proposition \ref{Prop estim comM2 s=1/2}
and we plug the estimate 
of Proposition \ref{com2M s=1/2}
in the differential inequation \eqref{ineq energ}.
We get the same estimate \eqref{dt Mf <}
and from Gronwall's lemma
we conclude the proof of Theorem~\ref{Th Sobolev}.
\qed

{\bf Proof of Corollary \ref{Corol Sobolev}.}

We consider the case $\frac12<s<1$.
We first note that, from Theorem \ref{Th Sobolev},
$f(t,\cdot) \in H^{+\infty}_2(\RR)$ for all $t>0$.
We introduce the mollifier
\begin{equation*}
M(\xi)=  1+\ \xi^2
\end{equation*}
which corresponds to the differential operator $M  =  1 - \Delta_v$.

By a proof similar to that of propositions \ref{Prop estim comM0}
and \ref{Prop estim comM2},
since $M$ satisfies obviously the estimates of lemma \ref{lem estime M},
we can prove the following estimates of the commutators:
For $f,g \in L^2_1$ and $h\in L^2(\RR)$, we have
\begin{align*}
\big| \big( M  K(f, g) , h \big)_{L^2}
  - \big( K(f, M g) , h \big)_{L^2} \big|
&\leq C \, \| f \|_{L^1} \|M g\|_{L^2} \|h\|_{L^2},
\\
\left|
\big(  (v^2 M) K(f, g) , h\big)_{L^2} 
- \big(K(f, (v^2 M) g) , h\big)_{L^2}
\right|
&\leq
  C \, \|f\|_{L^1_2} \,
  \|M g\|_{H^{\frac12}_2} \, 
  \|h\|_{H^{\frac12}}
\end{align*}
where $C$ depends only on $\beta$ and
$\| f \|_{  L^{\infty}( ]0,+\infty[ ; L^1_{2+2s} \bigcap L\log L(\RR) ) }$.

Following the proof of Theorem \ref{Th Sobolev}
and the same notations,
we get a differential inequation similar to \eqref{ineq energ}
(remark that the mollifier $M$ is independent of the time)
\begin{align*}
\frac12 \frac{d}{dt}
\Bigl( \| M f \|_{L^2}^2 + \| v^2 M f \|_{L^2}^2 \Bigr)
  +
c_f \| M f \|_{H^s_2}^2
\leq
\hbox{com}_0 + \hbox{com}_2
+ C \| f \|_{L^1} \| M f \|_{L^2_2}^2.
\end{align*}
The previous estimates of the commutators imply
\begin{align*}
\frac12 \frac{d}{dt}
\Big( 
  \|M f\|_{L^2}^2 + \|v^2 M  f\|_{L^2}^2 
 \Big)
  + c_f \|M f\|_{H^s_2}^2 
\leq
C \|f\|_{L^1_2} \|M f\|_{H^{\frac12}_2}^2.
\end{align*}
Using the interpolation estimate \eqref{interp H1/2}
we deduce the following differential inequation
\begin{equation*}
\frac12 \frac{d}{dt}
\Big( \|M f\|_{L^2}^2 + \|v^2 M f\|_{L^2}^2 \Big)
\leq
C\Big( \|M f\|_{L^2}^2 + \|v^2 M f\|_{L^2}^2 \Big)
\end{equation*}
and from Gronwall's lemma we derive
\begin{equation*}
\|M f\|_{L^2_2}^2
\le
 e^{2 C t} 
 \|M f_0\|_{L^2_2}^2.
\end{equation*}
That concludes
the proof of Corollary \ref{Corol Sobolev} in the case $\frac12<s<1$.

The proof in the case $s=\frac12$ is similar.
\qed

\renewcommand{\theequation}{\thesection.\arabic{equation}}
\setcounter{equation}{0}
\section{ANALYTICITY PROPERTY FOR KAC'S EQUATION}\label{section Kac}
From the Theorem \ref{Th Sobolev},
we know that the weak solution
of the Cauchy problem of the Kac's equation \eqref{eq Kac}
has the following regularity: for any $t_0>0$,
$f\in L^{\infty}([t_0,T_0];H^2_2(\RR))$.
Therefore $f$ is a solution of the following Cauchy problem :
\begin{equation*}
\left\{
\begin{array}{ll}
  & \displaystyle\frac{\partial f}{\partial t}
    = K(f,f) ,\\
  & f|_{t=0} = f_0 \in H^2_2(\RR).
\end{array}
\right.
\end{equation*}
and we can suppose that the initial datum is 
$f_0 \in H^2_2(\RR) \bigcap L^1_2(\RR)$.

We have the local analytic regularizing
effect of Cauchy problem.
\begin{theo}\label{Th analytic}
Assume that the cross-section kernel $\beta$ satisfies
\eqref{beta},
the initial datum 
$f_0 \in L^1_{2+2s} \bigcap H^2_2(\RR)$
and $f\in L^{\infty}([0,T_0]; H^2_2 \bigcap L^1_2(\RR))$
is a nonnegative weak solution 
of the Cauchy problem of the Kac's equation \eqref{eq Kac}
for some $T_0>0$.

\noindent
{\bf - Case} $\frac12 < s < 1$.

There exist $0<T_* \leq T_0$ and
$c_0>0$ such that
\begin{equation*}
e^{c_0 t<D_v>} f
\in L^{\infty}([0,T_*];L^2_1(\RR)).
\end{equation*}
Therefore we have 
$f(t,\cdot) \in G^1(\RR)$ for any $0<t\leq T_*$.

\noindent
{\bf - Case} $s=\frac12$.

For any $0<\alpha<1$,
there exist $0<T_* \leq T_0$ and
$c_0>0$ such that
\begin{equation*}
e^{c_0 t<D_v>^\alpha} f
\in L^{\infty}([0,T_*];L^2_1(\RR)).
\end{equation*}
Therefore for any $0<\alpha<1$ and $0<t\leq T_*$ we have
$f(t,\cdot) \in G^{1/\alpha}(\RR)$.
\end{theo}

{\bf Proof of the Theorem \ref{Th analytic}.}

We choose the test function 
\begin{equation*}
\tilde{\varphi}(t,\cdot)
= 
\big(G_{\delta}(t,D_v)
  \langle v \rangle^2G_{\delta}(t,D_v)f\big)(t,\cdot)
\in L^{\infty}(]0,T_0[;H^2(\RR))
\end{equation*}
where the mollifier $G_{\delta}$ 
is given in section \ref{section Estimate} by \eqref{G_delta}.

We have
\begin{equation*}
\left(\partial_t f \,,\,\tilde{\varphi} \right)_{L^2}
= 
\left(K(f,f)\,,\,\tilde{\varphi} \right)_{L^2} .
\end{equation*}
First, the left-hand side term is 
\begin{align*}
\left(
  \partial_t f \,,\, \tilde{\varphi}
\right) _{L^2} 
= \frac12 \frac{d}{dt}
  \| G_{\delta} f \|^2_{L^2_1}
- 
\left(
  (\partial_{t} G_{\delta}) f \,,\, 
  G_{\delta} f
\right)_{L^2}
-
\left(
  v(\partial_{t} G_{\delta})f) \,,\,
  v G_{\delta} f
\right)_{L^2}. 
\end{align*}
The rights-hand side is
\begin{align*}
\left(K(f,f)\,,\,\tilde{\varphi} \right)_{L^2} 
&= 
\left(
  G_{\delta}(K(f,f) \,,\,
  (1+ v^2)G_{\delta} f 
\right)_{L^2}
\\
&= 
\left(
  K(f,G_{\delta}f) \,,\, G_{\delta} f 
\right)_{L^2} 
+ 
\left(
    K(f,vG_{\delta}f) \,,\, vG_{\delta} f
\right)_{L^2}
\\
& + 
\left(
  G_{\delta} K(f,f)- K(f,G_{\delta}f) \,,\,
  G_{\delta} f 
\right)_{L^2}
\\
&+ 
\left( 
  v G_{\delta} K(f,f)- K(f,vG_{\delta}f) \,,\,
  v G_{\delta} f 
\right)_{L^2}.
\end{align*}
Therefore we can write :
\begin{align}\label{d/dtGf}
\frac12 \frac{d}{dt}
  \| G_{\delta} f \|^2_{L^2_1}
&- \left(
  K(f,G_{\delta}f) \,,\, G_{\delta} f 
\right)_{L^2} 
- 
\left(
    K(f,vG_{\delta}f) \,,\, vG_{\delta} f
\right)_{L^2}
\\
&=
\hbox{(time term)} + \hbox{(commutator)}\nonumber
\end{align}
where
\begin{align*}
\hbox{(time term)}
=
& - \left(
  (\partial_{t} G_{\delta})(\cdot,D_v) f \,,\, 
  G_{\delta}(t,D_v) f(t,\cdot)
\right)_{L^2}
\\
& - \left(
  v(\partial_{t} G_{\delta})(\cdot,D_v)f \,,\,
  v G_{\delta}(t,D_v)f
\right)_{L^2}
\end{align*}
and
\begin{align*}
\hbox{(commutator)}
= 
&\left(
  G_{\delta} K(f,f)- K(f,G_{\delta}f) \,,\,
  G_{\delta} f 
\right)_{L^2}
\\
&+ 
\left( 
  v G_{\delta} K(f,f)- K(f,vG_{\delta}f) \,,\,
  v G_{\delta} f 
\right)_{L^2}.
\end{align*}
Furthermore, 
\begin{equation*}
\|G_{\delta}f\|^2_{H^s_1} 
\leq
\|G_{\delta}f\|^2_{H^s} 
+ \|vG_{\delta}f\|^2_{H^s} 
+ \|G_{\delta}f\|^2_{L^2}
\end{equation*}
and by the Proposition \ref{Prop coer}, 
\begin{equation*}
- \left( K(f,G_{\delta}f) , (G_{\delta}f) \right)_{L^2} 
\geq
c_f \|G_{\delta}f\|^2_{H^{s}} - C\|f\|_{L^1}\|G_{\delta}f\|^2_{L^2} 
\end{equation*}
and
\begin{equation*}- \left( K(f,vG_{\delta}f) , v(G_{\delta}f) \right)_{L^2}
\geq
c_f\|vG_{\delta}f\|^2_{H^{s}} - 
C\|f\|_{L^1}\|vG_{\delta}f\|^2_{L^2}.
\end{equation*}
Then the equality \eqref{d/dtGf} became
\begin{align}\label{d/dtGf<}
\frac12 \frac{d}{dt}
  \|  G_{\delta} f \|^2_{L^2_1}
& + c_f\|G_{\delta}f\|^2_{H^{s}_1}
\\
&\leq
 C\|f\|_{L^1}\|G_{\delta}f\|^2_{L^2_1}
+
\hbox{(time term)} + \hbox{(commutator)}\nonumber.
\end{align}

\noindent
{\bf - Case :} $\frac12<s<1$.
\\
\noindent
We consider the mollifier $G_{\delta}$ 
defined in \eqref{G_delta}
and chosen with $\alpha=1$
\begin{equation*}
   G_\delta(t, \xi) =
  \frac{     e^{c_0 t \langle \xi \rangle } }
  {1+ \delta e^{c_0 t \langle \xi \rangle } }.
\end{equation*}

\begin{rem}
This is the optimal choice for $\alpha\in]0,2[$
as it can be seen in the estimates of section \ref{section Estimate} : 
for example, from lemma \ref{lem I2b},
the term
$
\|G_\delta f\|_{H^{(\frac{3\alpha}2-1)^+}}
\|G_\delta f\|_{H^{\frac{\alpha}2}}
$
can be controlled by
the coercivity
only if $\alpha\leq1$.
\end{rem}

Using the Propositions 
\ref{Prop estim com0}, \ref{Prop estim com1}
and the lemma \ref{lem temp}
of section \ref{section Estimate}, we get :
\\
{\sl - Estimate of commutators terms:}
\begin{equation}\label{com0-1}
| \left(
  G_\delta K(f, f)  
    - K(f, G_\delta f) 
  \,,\, 
  G_\delta f
  \right)_{L^2} |
\leq
  C \|G_\delta f\|_{L^2_1} 
    \|G_\delta f\|_{H^{1/2}}^2 
\end{equation}
and
\begin{equation}\label{com1-1}
|
  \left(
  v G_\delta) K(f, f)  
  - K(f, (v G_\delta) f)
  \,,\, 
  v G_\delta f
  \right)_{L^2} |
\leq
 C 
( \| f \|_{L^1_2} + \|G_\delta f\|_{L^2_1}) 
\|G_\delta f \|_{H^{\frac12}_1}^2  .
\end{equation}
\noindent
{\sl - Estimate of the terms involving the derivative
 with respect to time:}
\begin{equation}\label{estim time0-1}
\left|\left(
  (\partial_{t} G_{\delta})(t,D_v)f(t,\cdot) \,,\,
  G_{\delta}(t,D_v)f(t,\cdot)
\right)_{L^2} \right|
\leq 
  C \, \|G_{\delta}f\|^2_{H^{1/2}},
\end{equation}
and
\begin{equation}\label{estim time1-1}
\left|\left(
  v(\partial_{t} G_{\delta})(t,D_v)f(t,\cdot) \,,\,
  vG_{\delta}(t,D_v)f(t,\cdot)
\right)_{L^2} \right|
\leq 
  C \, 
\|G_{\delta}f \|^2_{H^{1/2}_1}.
\end{equation}
Therefore, plugging the estimates \eqref{com0-1}-\eqref{estim time1-1}
into \eqref{d/dtGf<},
we obtain
\begin{equation*}
\frac12 \frac{d}{dt}\|G_{\delta}f\|^2_{L^2_1} 
  + c_f \|G_{\delta}f\|^2_{H^s_1}
\leq
  C \| G_{\delta}f \|^2_{H^{\frac12}_1} 
+ C \| G_{\delta}f \|_{L^2_1} 
    \| G_{\delta}f \|^2_{H^{\frac12}_1}.
\end{equation*}
From the interpolation inequality \eqref{interp H1/2} we have
\begin{equation*}
C \| G_\delta f \|^2_{H^{\frac12}_1}
\leq 
C \lambda_1 \| G_\delta f \|^2_{H^s_1} 
+
C \lambda_1^{\frac{-1}{2s-1}} \| G_\delta f \|^2_{L^2_1}
\end{equation*}
and
\begin{equation*}
C
\| G_{\delta}f \|_{L^2_1} 
\| G_\delta f \|^2_{H^{\frac12}_1}
\leq
C \lambda_2 
\| G_{\delta}f \|_{L^2_1} 
\| G_\delta f \|^2_{H^s_1} 
+
C
\lambda_2^{\frac{-1}{2s-1}} 
\| G_\delta f \|^3_{L^2_1}.
\end{equation*}
Choosing $\lambda_1$ and $\lambda_2$ such that
$C \lambda_1 = c_f/4$
and
$C \lambda_2 \| G_{\delta}f \|_{L^2_1} = c_f/4$
we get
\begin{align*}
\frac12 \frac{d}{dt} \| G_{\delta} f \|^2_{L^2_1} 
&\leq
C_1 \| G_\delta f \|^2_{L^2_1} 
+ 
C_2 \| G_{\delta} f \|^{\frac{1}{2s-1}+3}_{L^2_1}
\end{align*}
for $t \in [0,T_0]$, with $C_1,\,C_2 >0$ 
are independent of $t$ and $\delta>0$.
Then
\begin{align*}
\frac{d}{dt} \| G_{\delta} f \|_{L^2_1} 
&\leq
C_1 \| G_\delta f \|_{L^2_1} 
+ 
C_2 \| G_{\delta} f \|^{\gamma}_{L^2_1}
\end{align*}
where
$\gamma = \frac{1}{2s-1}+2$.
We set $\psi(t) = \| G_{\delta} f(t,\cdot) \|_{L^2_1}$.
Therefore
\begin{equation*}
\frac{d}{dt} \psi(t) \leq C_1 \psi(t) + C_2 \psi(t)^\gamma.
\end{equation*}
Solving the previous differential inequation,
we easily get
\begin{equation*}
\psi(t) 
\leq
\frac
{ e^{C_1 t} \, \psi(0)}   
{
\left(
      1 - \frac{C_2}{C_1} 
      \left( e^{(\gamma-1) C_1 t} -1 \right)
      \psi(0)^{\gamma-1}
\right)^{\frac1{\gamma-1}}
},
\end{equation*}
that is
\begin{equation*}
 \| G_\delta f(t,\cdot) \|_{L^2_1}
\leq
\frac
{ e^{C_1 \, t} \, \| f_0 \|_{L^2_1} }   
{
\left(
      1 - \frac{C_2}{C_1} 
      \left( e^{\frac{s}{s-\frac12} \, t} -1 \right)
      \| f_0 \|_{L^2_1}^{\frac{s}{s-\frac12}}
\right)^{\frac{s-\frac12}{s}}
}.
\end{equation*}
We now choose $0<T_* \leq T_0$ small enough
so that for $t\in[0,T_*]$
\begin{equation}\label{T_*}
1 - \frac{C_2}{C_1} 
\left( e^{\frac{s}{s-\frac12} \, t} -1 \right)
\| f_0 \|_{L^2_1}^{\frac{s}{s-\frac12}}
\geq \left(\frac12\right)^{\frac{s}{s-\frac12}},
\end{equation}
and taking $\delta \rightarrow 0$, we have for $t\in]0,T_*]$,
\begin{equation*}
\|
  e^{c_{0} t <D_v>} f(t,\cdot)
\|_{L^2_1} 
\leq 
2 e^{C_1 t} \|f_0\|_{L^2_1}.
\end{equation*}
This concludes the proof of Theorem \ref{Th analytic}
in the case $\frac12<s<1$.

\noindent
{\bf - Case :} $s=\frac12$.
\\
\noindent
We consider the mollifier $G_{\delta}$ defined in \eqref{G_delta}
with $0<\alpha<1$.
Taking ${\alpha'}=\frac12$ in the estimate of the commutator
in Proposition \ref{Prop estim com1 s=1/2}
we obtain
\begin{align*}
\frac12 \frac{d}{dt}\|G_{\delta}f\|^2_{L^2_1} 
  &+ c_f \|G_{\delta}f\|^2_{H^{\frac12}_1}
\leq
\\
  &C \| f \|_{L^1_2} 
    \| G_{\delta}f \|^2_{H^{\frac14}_1} 
+ C \| G_{\delta}f \|_{L^2_1} 
    \| G_{\delta}f \|^2_{H^{\frac{\alpha}2}_1}
+C  \| G_{\delta}f \|^2_{H^{\frac{\alpha}2}_1}.
\end{align*}
From an interpolation estimate similar as \eqref{interp H1/2},
we get the following differential inequation
\begin{align*}
\frac12 \frac{d}{dt}\|G_{\delta}f\|^2_{L^2_1} 
\leq
C_1' \| G_{\delta} f \|_{L^2_1}^2
+ C_2' \| G_{\delta}f \|_{L^2_1}^{\frac{\alpha}{1-\alpha}+3}.
\end{align*}
where $C_1' , C_2' > 0$ are independent of $\delta > 0$.
This concludes the proof of the Theorem~\ref{Th analytic}.
\qed

{\bf Proof of the propagation of analyticity and end of the proof of Theorem \ref{Th Kac}.}

We could use the Theorem 2.6 of \cite{DFT09}
(propagation of Gevrey regularity in the case of an even initial datum $f_0$).
We present below a direct proof.

We consider the case $\frac12 <s < 1$.
Let $f$ a nonnegative weak solution of the Cauchy problem \eqref{eq Kac}
which fulfills the assumptions of Theorem \ref{Th Kac}:
From Theorem~\ref{Th Sobolev}, we have $f(t,\cdot)\in H^2_2(\RR)$
for all $t>0$.

Let us consider some arbitrary and fixed $0<T_0<T_1$.
From Corollary \ref{Corol Sobolev},
there exists a constants $C_0$ which depends only on 
$\beta$, 
$\| f \|_{  L^{\infty}( ]0,+\infty[ ; L^1_{2+2s} \bigcap L\log L(\RR) ) }$,
$T_0$, $T_1$ and $\| f(T_0,\cdot) \|_{H^2_2}$ such that
\begin{equation}\label{estim H^2_2}
\forall t \in[T_0, \, T_1+1], \quad
\| f(t,\cdot) \|_{H^2_2}
\leq
C_0.
\end{equation}
We then consider $t_0\in [T_0, T_1]$ and we apply Theorem \ref{Th analytic}
for the initial value $\tilde f_0 = f(t_0, \cdot)$ and for $t\in[0, 1]$.
There exist $0<T_*\leq 1$ and some constants $c_0>0$ and $C_1>0$  such that
\begin{equation*}
\forall t\in ]t_0, t_0+T_*], \quad
 \| e^{c_{0} (t-t_0) <D_v>} f(t,\cdot)
\|_{L^2_1} 
\leq 
2 e^{C_1 (t-t_0)} \| f(t_0,\cdot) \|_{L^2_1}.
\end{equation*}
In addition, we remark from the proof of Theorem \ref{Th analytic}
and the inequality \eqref{T_*}
that the time $T_*$ depends only on $T_0$, $T_1$
and $\| f \|_{  L^{\infty}( [T_0,T_1+1] ; H^2_2 \bigcap L^1_2(\RR) ) }$,
which is controlled 
by $\| f_0 \|_{L^1_2}$ and by the constant $C_0$ from \eqref{estim H^2_2},

Hence $T_*>0$ can be chosen independent of $t_0\in [T_0, T_1]$.
Therefore $f(t,\cdot) \in G^1(\RR)$ for all $t\in]t_0, t_0+T_*]$,
and this will remain true for $t \in ]T_0, T_1]$.
Since $T_0<T_1$ are arbitrary, 
the solution of the Cauchy problem \eqref{eq Kac} satisfies
$f(t,\cdot) \in G^1(\RR)$ for any $t>0$.

The proof for the case $s=\frac12$ is similar.
This concludes the proof of Theorem \ref{Th Kac}.
\qed

\renewcommand{\theequation}{\thesection.\arabic{equation}}
\setcounter{equation}{0}
\section{ANALYTICITY PROPERTY FOR BOLTZMANN EQUATION}\label{section Boltzmann}

In this section, we will prove the analyticity 
of the radially symmetric solutions
of the Boltzmann equation (Theorem \ref{Th Boltzmann})
\begin{equation*}
\frac{\partial f}{\partial t} = Q(f,f),\quad
v\in \RR^{3},\, t>0; \quad f|_{t= 0} = f_0.
\end{equation*}
Using the Bobylev's formula, we have for $\xi\in\RR^3$
\begin{equation*}
\cF\left(Q(f,g) \right) (\xi)
= 
\int_{S^2}
b\left(\frac{\xi}{|\xi|} \,.\, \sigma\right)
\left\lbrace 
  \hat{f}(\xi^-)\hat{g}(\xi^+)-\hat{f}(0)\hat{g}(\xi)
\right\rbrace 
d\sigma
\end{equation*}
where
\begin{equation*}
\xi^+ = \frac{\xi+|\xi| \, \sigma}2, \quad
\xi^- = \frac{\xi-|\xi| \, \sigma}2.
\end{equation*}
We define $\theta$ by
\begin{equation*}
  \cos\theta = \langle \frac{\xi}{|\xi|} \,,\, \sigma \rangle.
\end{equation*}
We then have
\begin{equation*}
|\xi^+| = |\xi| \, |\cos\textstyle{\frac{\theta}2|}, \quad
|\xi^-| = |\xi| \, |\sin\textstyle{\frac{\theta}2}|.
\end{equation*}
Let $f$ be a radially symmetric function.
That is $f(A v) = f(v)$ for any orthogonal $3 \times 3$
matrix $A$.
Therefore $f(v) = f(0,0,|v|)$.
We compute the Fourier transform in $\RR^3$ :
\begin{align*}
\cF_{\RR^3}(f)(\xi)
  = \int_{\RR^3} e^{-i \xi.v} f(v) dv
  = \int_{\RR^3} e^{-i |\xi| v_3} f(v) dv
  = \int_{\RR}  e^{-i |\xi| u} F(u) du
\end{align*}
where
\begin{equation}\label{F=int f}
F(u) = \int_{\RR^2}  f(v_1, v_2, u) dv_1 dv_2.
\end{equation}
Then
\begin{equation}\label{Ff=FF}
\hat f (\xi) = \hat F(|\xi|)
\end{equation}
is also radially symmetric.
Let us consider two radially symmetric functions $f$ and $g$
and let us denote $F$ and $G$ the associated functions
defined by \eqref{F=int f}.
We compute
\begin{align*}
\cF\left(Q(f,g)\right) (\xi)
&=
\int_{S^2}
b\left(\frac{\xi}{|\xi|} \,.\, \sigma\right)
\left\lbrace
  \hat{F}(|\xi^-|)\hat{G}(|\xi^+|)-\hat{F}(0)\hat{G}(|\xi|)
\right\rbrace
d\sigma
\\
&=
\int_{-\frac{\pi}2}^{\frac{\pi}2}
\beta(\theta)
\left\lbrace
  \hat{F}(|\xi^-|)\hat{G}(|\xi^+|)-\hat{F}(0)\hat{G}(|\xi|)
\right\rbrace
d\theta
\end{align*}
where
\begin{equation*}
\beta(\theta) = 2 \pi \sin\theta \, b(\cos\theta).
\end{equation*}
Let $f(t,\cdot)$ be a solution of the Boltzmann equation.
We put for $t\geq0$
\begin{equation*}
F(t,u) = \int_{\RR^2} f(t, v_1, v_2, u) dv_1 dv_2.
\end{equation*}
Therefore $\hat f(t,\cdot)$ is solution of the equation
\begin{equation*}
 \partial_t \hat f(t,\xi)
  = \cF \big( Q(f(t,\cdot), f(t,\cdot)) \big) (\xi)
\end{equation*}
and from \eqref{Ff=FF}
we prove that $F(t,\cdot)$ is a solution of the equation
\begin{equation*}
 \partial_t \hat F(t,|\xi|)
  =
\int_{-\frac{\pi}2}^{\frac{\pi}2}
\beta(\theta)
\left\lbrace
  \hat{F}(t,|\xi\sin\textstyle{\frac{\theta}2}| )
  \hat{F}(t,|\xi\cos\textstyle{\frac{\theta}2}| )
 -\hat{F}(t,0)\hat{F}(t,|\xi|)
\right\rbrace
d\theta.
\end{equation*}
We use the following lemma (see \cite{LX09}):
\begin{lem}\label{lem reg dim3 dim1}
Let $f\in L^1_k(\RR^3)$ radially symmetric, $f\geq 0$
and define $F$ by
\eqref{F=int f}.
Then $F\in L^1_k(\RR)$.
Assume that $f$ is also uniformly integrable $f\geq 0$.
Then $F$ is also uniformly integrable.
\end{lem}

From lemma \ref{lem reg dim3 dim1},
$F(t, \cdot)\in L^1_{2+2s}(\RR)$, but
we do not have $F(t, \cdot)\in L\log L(\RR)$.
However $F$ is uniformly integrable, and it is enough
to get the coercivity property of propo\-sition~\ref{Prop coer}.

{\bf Proof of Theorem \ref{Th Boltzmann}}.

\noindent
Case $\frac12<s<1$ (the proof for the case $s=\frac12$ is similar).
We apply the Theorem \ref{Th Kac} : 
for a fixed $t>0$, there exists a constant $c_0$ such that
\begin{equation*}
\| e^{c_0 \langle \cdot \rangle} \hat F(t,|\cdot|) \|_{L^2(\RR)}
< \infty
\end{equation*}
that is
$F(t,\cdot) \in G^1(\RR)$.
Since 
\begin{equation*}
e^{c_0 \langle \xi \rangle} \hat F(t,|\xi|) 
= 
e^{c_0 \langle \xi \rangle} \hat f(t,\xi) 
\end{equation*}
we have
\begin{equation*}
\| e^{\frac{c_0}2 \langle \cdot \rangle } \hat f(t,\cdot) \|_{L^2(\RR^3)}
\leq
C \, \| e^{c_0 \langle \cdot \rangle } \hat f(t, |\cdot|) \|_{L^2(\RR)}
 < \infty
\end{equation*}
and therefore
$f(t,\cdot) \in G^1(\RR^3)$.

The proof for the case $s=\frac12$ is similar.
This concludes the proof of Theorem \ref{Th Boltzmann}.
\qed

\medskip

\noindent
{\bf Acknowledgements :}
The authors wish to thank Chao-Jiang Xu for interesting conversations.

\end{document}